\documentclass[12pt]{amsart}
\usepackage{amsfonts}
\usepackage{amsmath,amssymb,amsthm}

\newcommand{\R}{\mathbb R}

\newcommand{\E}{\mathbb E}

\renewcommand{\span}{\mathrm{span}}
\newcommand{\tr}{\mathrm{tr}}

\newtheorem{thm}{Theorem}[section]

\newtheorem{lem}[thm]{Lemma}
\newtheorem{prop}[thm]{Proposition}
\theoremstyle{definition}
\newtheorem{defn}[thm]{Definition}
\theoremstyle{remark}

\newcommand{\ds}{\displaystyle}

\begin{document}

\title[THEORY OF SPACELIKE SURFACES IN MINKOWSKI SPACE]
{AN INVARIANT THEORY OF SPACELIKE SURFACES IN THE FOUR-DIMENSIONAL
MINKOWSKI SPACE}

\author{Georgi Ganchev and Velichka Milousheva}
\address{Bulgarian Academy of Sciences, Institute of Mathematics and Informatics,
Acad. G. Bonchev Str. bl. 8, 1113 Sofia, Bulgaria}
\email{ganchev@math.bas.bg}
\address{Bulgarian Academy of Sciences, Institute of Mathematics and Informatics,
Acad. G. Bonchev Str. bl. 8, 1113, Sofia, Bulgaria; "L. Karavelov"
Civil Engineering Higher School, 175 Suhodolska Str., 1373 Sofia,
Bulgaria} \email{vmil@math.bas.bg}


\begin{abstract}
We consider spacelike surfaces in the four-dimensional Minkowski
space and introduce geometrically an invariant linear map of
Weingarten-type in the tangent plane at any point of the surface
under consideration. This allows us to introduce principal lines
and an invariant moving frame field. Writing derivative formulas
of Frenet-type for this frame field, we obtain eight invariant
functions. We prove a fundamental theorem of Bonnet-type, stating
that these eight invariants under some natural conditions
determine the surface up to a motion.

We show that the basic geometric classes of spacelike surfaces in
the four-dimensional Minkowski space, determined by conditions on
their invariants, can be interpreted in terms of the properties of
the two geometric figures: the tangent indicatrix, and the normal
curvature ellipse.

We apply our theory to a class  of spacelike general rotational
surfaces.
\end{abstract}

\maketitle

\section{Introduction}

In this paper we consider the general theory of spacelike surfaces
in the four-dimensional Minkowski space $\R^4_1$. The basic
feature of our approach to this theory is the introduction of an
invariant linear map of Weingarten-type in the tangent plane at
any point of the surface. Studying surfaces in the Euclidean space
$\R^4$, in \cite{GM1} we introduced a linear map of
Weingarten-type, which plays a similar role in the theory of
surfaces in $\R^4$ as the Weingarten map in the theory of surfaces
in $\R^3$. We gave a geometric interpretation of the second
fundamental form and the Weingarten map of the surface in
\cite{GM4}. Following our approach to the surfaces in $\R^4$, here
we develop the theory of spacelike surfaces in $\R^4_1$ in a
similar way.

Let $M^2$  be a spacelike surface in $\R^4_1$. Considering the
tangent space $T_pM^2$ at a point $p \in M^2$, we introduce an
invariant $\zeta_{\,g_1,g_2}$ of a pair of two tangents $g_1$,
$g_2$ using the second fundamental tensor $\sigma$ of $M^2$. By
means of this invariant we define conjugate, asymptotic, and
principal tangents.

The second fundamental form $II$ of the surface $M^2$ at a point
$p \in M^2$ is introduced on the base of conjugacy of two tangents
at the point. The second fundamental form $II$ determines an
invariant linear map of Weingarten-type $\gamma: T_pM^2
\rightarrow T_pM^2$ at any point of $M^2$ in the standard way. The
map $\gamma$ generates two invariants  $k$ and $\varkappa$. We
prove that the invariant $\varkappa$ is the curvature of the
normal connection of the surface. The number of asymptotic
tangents at a point of $M^2$ is determined by the sign of the
invariant $k$. In the case $k=0$ there exists a one-parameter
family of asymptotic lines, which are principal. It is interesting
to note that the "umbilical" points, i.e. points at which the
coefficients of the first and the second fundamental forms are
proportional, are exactly the points at which the mean curvature
vector $H$ is zero. Minimal spacelike surfaces are characterized
by the equality $\varkappa^2 - k =0$.

Analogously to $\R^3$, the invariants $k$ and $\varkappa$ divide
the points of $M^2$ into four types: flat, elliptic, hyperbolic
and parabolic points. The surfaces consisting of flat points  are
characterized by the conditions $k= \varkappa = 0$. In Section
\ref{S:flat points} we give a local geometric description of
spacelike surfaces consisting of flat points whose mean curvature
vector at any point is a non-zero spacelike vector or  timelike
vector. We prove that:

\emph{Any spacelike surface consisting of flat points whose mean
curvature vector at any point is a non-zero spacelike vector or
timelike vector  either lies in a hyperplane of $\R^4_1$ or is
part of a developable ruled surface in $\R^4_1$}.

We introduce the indicatrix of Dupin $\chi$ at an arbitrary
(non-flat) point of a spacelike surface in $\R^4_1$ by means of
the second fundamental form as in the theory of surfaces in
$\R^3$. Then the elliptic, hyperbolic and parabolic points of a
spacelike surface $M^2$ in $\R^4_1$ are characterized in terms of
the indicatrix $\chi$ as in $\R^3$. The conjugacy in terms of the
second fundamental form coincides with the conjugacy with respect
to the indicatrix $\chi$. In Section \ref{S:Weingarten} we prove
that:

\emph{A spacelike surface $M^2$ is minimal if and only if the
indicatrix $\chi$ is a circle.}

\emph{A spacelike surface $M^2$ is with flat normal connection if
and only if the indicatrix $\chi$ is a rectangular hyperbola (a
Lorentz circle).}

In the local theory of surfaces in Euclidean space  a statement of
significant importance is a theorem of Bonnet-type giving the
natural conditions under which the surface is determined up to a
motion. A theorem of this type was proved for surfaces with flat
normal connection by B.-Y. Chen in \cite{Chen1}. In \cite{GM2} we
proved a fundamental theorem of Bonnet-type for minimal surfaces
in $\R^4$ and in \cite{GM3} we proved such a theorem for surfaces
in $\R^4$ free of minimal points. Here we consider spacelike
surfaces in $\R^4_1$ whose mean curvature vector at any point is a
non-zero spacelike vector or  timelike vector. Using a
geometrically determined moving frame of Frenet-type on such a
surface  and the corresponding derivative formulas, we obtain
eight invariant functions. In Section \ref{S:spacelikeH} and
Section \ref{S:timelikeH} we prove our basic Theorem \ref{T:Main
Theorem} and Theorem  \ref{T:Main Theorem-2}, stating that

\emph{Any spacelike surface whose mean curvature vector at any
point is a non-zero spacelike vector or  timelike vector is
determined up to a motion in $\R^4_1$ by its eight invariant
functions satisfying some natural conditions.}

In Section \ref{S:Examples} we apply our theory to a class of
spacelike general rotational surfaces in $\R^4_1$.

\section{Invariants of a tangent on a spacelike surface in $\R^4_1$}
\label{S:Invariants}

We consider the Minkowski space $\R^4_1$ endowed with the metric
$\langle , \rangle$ of signature $(3,1)$. Let $Oe_1e_2e_3e_4$ be a
fixed orthonormal coordinate system in $\R^4_1$, i.e. $e_1^2 =
e_2^2 = e_3^2 = 1, \, e_4^2 = -1$, giving the orientation of
$\R^4_1$.

A surface $M^2: z = z(u,v), \, \, (u,v) \in {\mathcal D}$
(${\mathcal D} \subset \R^2$) in $\R^4_1$ is said to be
\emph{spacelike} if $\langle , \rangle$ induces  a Riemannian
metric $g$ on $M^2$. Thus at each point $p$ of a spacelike surface
$M^2$ we have the following decomposition:
$$\R^4_1 = T_pM^2 \oplus N_pM^2$$
with the property that the restriction of the metric $\langle ,
\rangle$ onto the tangent space $T_pM^2$ is of signature $(2,0)$,
and the restriction of the metric $\langle , \rangle$ onto the
normal space $N_pM^2$ is of signature $(1,1)$.

Let $M^2$ be a spacelike surface in $\R^4_1$. The tangent space to
$M^2$ at an arbitrary point $p=z(u,v)$ of $M^2$ is ${\rm span}
\{z_u, z_v\}$, where $\langle z_u,z_u \rangle > 0$, $\langle
z_v,z_v \rangle > 0$. We use the standard denotations
\;$E(u,v)=\langle z_u,z_u \rangle, \; F(u,v)=\langle z_u,z_v
\rangle, \; G(u,v)=\langle z_v,z_v \rangle$ for the coefficients
of the first fundamental form
$$I(\lambda,\mu):= E \lambda^2 + 2F \lambda \mu + G \mu^2,\quad
\lambda, \mu \in \R.$$ Since $I(\lambda, \mu)$ is positive
definite we set $W=\sqrt{EG-F^2}$.

We choose a normal frame field $\{n_1, n_2\}$ such that $\langle
n_1, n_1 \rangle =1$, $\langle n_2, n_2 \rangle = -1$, and the
quadruple $\{z_u,z_v, n_1, n_2\}$ is positively oriented in
$\R^4_1$. Denote by $\nabla'$ the standard covariant derivative in
$\R^4_1$ and consider the functions
$$\begin{array}{ll}
\vspace{2mm}
c_{11}^1 = \langle z_{uu}, n_1 \rangle; & \qquad c_{11}^2 = \langle z_{uu}, n_2 \rangle;\\
\vspace{2mm}
c_{12}^1 = \langle z_{uv}, n_1 \rangle; & \qquad c_{12}^2 = \langle z_{uv}, n_2 \rangle;\\
\vspace{2mm} c_{22}^1 = \langle z_{vv}, n_1 \rangle; & \qquad
c_{22}^2 = \langle z_{vv}, n_2 \rangle.
\end{array} $$
Then we have the following derivative formulas:
$$\begin{array}{l}
\vspace{2mm} \nabla'_{z_u}z_u=z_{uu} = \Gamma_{11}^1 \, z_u +
\Gamma_{11}^2 \, z_v + c_{11}^1\, n_1 - c_{11}^2\, n_2;\\
\vspace{2mm} \nabla'_{z_u}z_v=z_{uv} = \Gamma_{12}^1 \, z_u +
\Gamma_{12}^2 \, z_v + c_{12}^1\, n_1 - c_{12}^2\, n_2;\\
\vspace{2mm} \nabla'_{z_v}z_v=z_{vv} = \Gamma_{22}^1 \, z_u +
\Gamma_{22}^2 \, z_v + c_{22}^1\, n_1 - c_{22}^2\, n_2,\\
\end{array}$$
where $\Gamma_{ij}^k$ are the Christoffel's symbols.
 If $\sigma$
denotes the second fundamental tensor of $M^2$, then we have
$$\begin{array}{l}
\sigma(z_u,z_u)=c_{11}^1\, n_1 - c_{11}^2\, n_2,\\
[2mm]
\sigma(z_u,z_v)=c_{12}^1\, n_1 - c_{12}^2\, n_2,\\
[2mm] \sigma(z_v,z_v)=c_{22}^1\, n_1 - c_{22}^2\, n_2.\end{array}
\leqno{(2.1)}$$

Obviously, the surface $M^2$ lies in a 2-plane if and only if
$M^2$ is totally geodesic, i.e. $c_{ij}^k=0, \; i,j,k = 1, 2.$ So,
we assume that at least one of the coefficients $c_{ij}^k$ is not
zero.

\vskip 2mm We shall define conjugate tangents at any point of the
surface $M^2$.

Let $g$ be a tangent at the point $p \in M^2$ determined by the
non-zero vector $X =  \lambda z_u + \mu z_v$. We consider the  map
$\sigma_g: T_pM^2 \rightarrow N_pM^2$, defined by
$$\sigma_g(Y) = \ds{\sigma \left(\frac{\lambda z_u + \mu z_v}{\sqrt{I(\lambda, \mu)}}, \,Y \right)}, \qquad Y \in T_pM^2. \leqno{(2.2)}$$
Obviously $\sigma_g$ is a linear map, which does not depend on the
choice of the non-zero vector $X$ collinear with $g$. Using (2.1)
and (2.2) we obtain the following decomposition of the normal
vectors $\sigma_g(z_u)$ and $\sigma_g(z_v)$:
$$\begin{array}{l}
\vspace{2mm}
\sigma_g(z_u)=\ds{\frac{\lambda\, c_{11}^1 + \mu \,c_{12}^1}{\sqrt{I(\lambda, \mu)}}\, n_1 - \frac{\lambda\, c_{11}^2 + \mu \,c_{12}^2}{\sqrt{I(\lambda, \mu)}}\, n_2},\\
\vspace{2mm} \sigma_g(z_v)=\ds{\frac{\lambda\, c_{12}^1 + \mu
\,c_{22}^1}{\sqrt{I(\lambda, \mu)}}\, n_1 - \frac{\lambda\,
c_{12}^2 + \mu \,c_{22}^2}{\sqrt{I(\lambda, \mu)}}\, n_2}.
\end{array} \leqno{(2.3)}$$

Let $g_1: X_1 =  \lambda_1 z_u + \mu_1 z_v$ and $g_2: X_2 =
\lambda_2 z_u + \mu_2 z_v$ be two tangents at the point $p \in
M^2$. We consider the parallelograms determined by the pairs of
normal vectors $\sigma_{g_1} (z_u)$, $\sigma_{g_2} (z_v)$ and
$\sigma_{g_2} (z_u)$, $\sigma_{g_1} (z_v)$  in the Lorentz plane
$\span\{n_1,n_2\}$. The oriented areas of these parallelograms are
denoted by $S(\sigma_{g_1} (z_u), \sigma_{g_2} (z_v))$, and
$S(\sigma_{g_2} (z_u), \sigma_{g_1} (z_v))$, respectively. We
assign the quantity $\zeta_{\,g_1,g_2}$ to the pair of tangents
$g_1$, $g_2$, defined by
$$\zeta_{\,g_1,g_2} = \ds{\frac{S(\sigma_{g_1} (z_u),\sigma_{g_2} (z_v))}{W} + \frac{S(\sigma_{g_2} (z_u), \sigma_{g_1} (z_v))}{W}}. \leqno{(2.4)}$$

Using equalities (2.3) we calculate that
$$\zeta_{\,g_1,g_2} = \ds{\frac{2\left|%
\begin{array}{cc}
\vspace{2mm}
  c_{11}^1 & c_{12}^1 \\
  c_{11}^2 & c_{12}^2 \\
\end{array}%
\right| \lambda_1 \lambda_2 + \left|%
\begin{array}{cc}
\vspace{2mm}
  c_{11}^1 & c_{22}^1 \\
  c_{11}^2 & c_{22}^2 \\
\end{array}%
\right| (\lambda_1 \mu_2 + \mu_1 \lambda_2) + 2\left|%
\begin{array}{cc}
\vspace{2mm}
  c_{12}^1 & c_{22}^1 \\
  c_{12}^2 & c_{22}^2 \\
\end{array}%
\right| \mu_1 \mu_2}{W \sqrt{I(\lambda_1, \mu_1)} \sqrt{
I(\lambda_2, \mu_2)}}}.$$

We introduce the following functions:
$$\Delta_1 = \left|%
\begin{array}{cc}
\vspace{2mm}
  c_{11}^1 & c_{12}^1 \\
  c_{11}^2 & c_{12}^2 \\
\end{array}%
\right|; \quad
\Delta_2 = \left|%
\begin{array}{cc}
\vspace{2mm}
  c_{11}^1 & c_{22}^1 \\
  c_{11}^2 & c_{22}^2 \\
\end{array}%
\right|; \quad
\Delta_3 = \left|%
\begin{array}{cc}
\vspace{2mm}
  c_{12}^1 & c_{22}^1 \\
  c_{12}^2 & c_{22}^2 \\
\end{array}%
\right|;$$
$$L(u,v) = \displaystyle{\frac{2 \Delta_1}{W}, \quad M(u,v) =
\frac{\Delta_2}{W}, \quad N(u,v) = \frac{2 \Delta_3}{W}}.$$ Hence,
$\zeta_{\,g_1,g_2}$ is expressed as follows:
$$\zeta_{\,g_1,g_2} = \ds{\frac{L \lambda_1 \lambda_2 + M (\lambda_1 \mu_2 + \mu_1 \lambda_2) + N \mu_1 \mu_2}{\sqrt{I(\lambda_1, \mu_1)}
\sqrt{I(\lambda_2, \mu_2)}}}. \leqno{(2.5)}$$

\vskip 3mm

\begin{prop}
The quantity $\zeta_{\,g_1,g_2}$ is invariant under any change of
the parameters on $M^2$.
\end{prop}

\vskip 2mm \noindent \emph{Proof}: Let
$$\begin{array}{l}
\vspace{2mm}
u = u(\bar u,\bar v);\\
\vspace{2mm} v = v(\bar u,\bar v),
\end{array}
\quad (\bar u,\bar v) \in \bar{\mathcal D}, \,\, \bar{\mathcal D}
\subset \R^2 \leqno{(2.6)}$$ be a smooth change of the parameters
$(u,v)$ on $M^2$ with $J = u_{\bar u} \, v_{\bar v} - u_{\bar v}
\, v_{\bar u}\neq 0$. Then
$$\begin{array}{l}
\vspace{2mm}
z_{\bar u} = z_u \,u_{\bar u} + z_v \,v_{\bar u},\\
\vspace{2mm} z_{\bar v} = z_u \,u_{\bar v} + z_v \,v_{\bar v}.
\end{array}$$
If $\bar E = \langle z_{\bar u}, z_{\bar u} \rangle$, $\bar F =
\langle z_{\bar u}, z_{\bar v} \rangle$ and $\bar G = \langle
z_{\bar v}, z_{\bar v} \rangle$, then we have
$$\begin{array}{l}
\vspace{2mm}
\bar E=u_{\bar u}^2\,E+2\,u_{\bar u}v_{\bar u}\,F+v_{\bar u}^2\,G,\\
\vspace{2mm} \bar F=u_{\bar u}u_{\bar v}\,E+(u_{\bar u}v_{\bar
v}+v_{\bar u}u_{\bar v})\,F
+v_{\bar u}v_{\bar v}\,G,\\
\vspace{2mm} \bar G=u_{\bar v}^2\,E+2\,u_{\bar v}v_{\bar
v}\,F+v_{\bar v}^2\,G\end{array}$$ and $\bar E \bar G - \bar
F^2=J^2\,(EG-F^2)$, hence $\bar W = \varepsilon J\,W, \,\,
\varepsilon = {\rm sign} \, J$.

Let
$$\begin{array}{l}
\vspace{2mm}
\sigma(z_{\bar u},z_{\bar u}) = \bar c_{11}^1 \, n_1 - \bar c_{11}^2 \, n_2,\\
\vspace{2mm}
\sigma(z_{\bar u},z_{\bar v}) = \bar c_{12}^1 \, n_1 - \bar c_{12}^2 \, n_2,\\
\vspace{2mm} \sigma(z_{\bar v},z_{\bar v}) = \bar c_{22}^1 \, n_1
- \bar c_{22}^2 \, n_2.
\end{array}$$
Then from (2.6) and (2.1) we find
$$\begin{array}{l}
\vspace{2mm} \bar c_{11}^k = u_{\bar u}^2 \,c_{11}^k  +
2 u_{\bar u}\, v_{\bar u}\, c_{12}^k + v_{\bar u}^2\,c_{22}^k,\\
\vspace{2mm} \bar c_{12}^k = u_{\bar u}\,u_{\bar v} \,c_{11}^k  +
(u_{\bar u} \, v_{\bar v}+ u_{\bar v}\, v_{\bar u})\, c_{12}^k +
v_{\bar u} \, v_{\bar v} \, c_{22}^k,\\
\vspace{2mm} \bar c_{22}^k = u_{\bar v}^2\,c_{11}^k  + 2 u_{\bar
v} \, v_{\bar v}\, c_{12}^k + v_{\bar v}^2\,c_{22}^k.
\end{array} \quad  \quad  (k = 1,2),$$
and hence
$$\begin{array}{l}
\vspace{2mm} \bar{\Delta}_1 = J \left(u_{\bar u}^2\,\Delta_1 +
u_{\bar u} \, v_{\bar u}\, \Delta_2 + v_{\bar u}^2\,\Delta_3 \right);\\
\vspace{2mm} \bar{\Delta}_2 = J \left(2 u_{\bar u}\,u_{\bar
v}\,\Delta_1 + (u_{\bar u}\,v_{\bar v} + u_{\bar v} \, v_{\bar
u})\,\Delta_2
+ 2 v_{\bar u}\, v_{\bar v}\,\Delta_3 \right);\\
\vspace{2mm} \bar{\Delta}_3 = J \left(u_{\bar v}^2\,\Delta_1 +
u_{\bar v} \, v_{\bar v}\, \Delta_2 + v_{\bar v}^2\,\Delta_3
\right).
\end{array}$$
Thus we find that the functions $\bar L$, $\bar M$, $\bar N$ are
expressed as follows:
$$\begin{array}{l}
\vspace{2mm}
\bar L=\varepsilon (u_{\bar u}^2\,L+2\,u_{\bar u}v_{\bar u}\,M+v_{\bar u}^2\,N),\\
\vspace{2mm} \bar M = \varepsilon (u_{\bar u}u_{\bar
v}\,L+(u_{\bar u}v_{\bar v}+v_{\bar u}u_{\bar v})\,M
+v_{\bar u}v_{\bar v}\,N),\\
\vspace{2mm} \bar N=\varepsilon (u_{\bar v}^2\,L+2\,u_{\bar
v}v_{\bar v}\,M+v_{\bar v}^2\,N).
\end{array}\leqno{(2.7)}$$
Hence, the functions $L,M,N$ change in the same way as the
coefficients of the first fundamental form $E,F,G$ under any
change of the parameters on $M^2$.

If $X =  \lambda z_u + \mu z_v = \bar \lambda z_{\bar u} + \bar
\mu z_{\bar v}$, then $\lambda = u_{\bar u} \bar \lambda + u_{\bar
v} \bar \mu, \,\,\mu = v_{\bar u} \bar \lambda + v_{\bar v} \bar
\mu$. Using (2.7) we obtain
$$\bar L \bar \lambda_1 \bar \lambda_2 +\bar M (\bar \lambda_1 \bar \mu_2 + \bar \mu_1 \bar \lambda_2) + \bar N \bar \mu_1 \bar \mu_2 =
\varepsilon \left(L \lambda_1 \lambda_2 + M (\lambda_1 \mu_2 +
\mu_1 \lambda_2) + N \mu_1 \mu_2 \right).$$ Having in mind that
$\bar I (\bar \lambda, \bar \mu) = I(\lambda, \mu)$, we get
$$\bar \zeta_{\,g_1,g_2} = \ds{\frac{\bar L \bar \lambda_1 \bar \lambda_2 + \bar M (\bar \lambda_1 \bar \mu_2 +
\bar \mu_1 \bar \lambda_2) + \bar N \bar \mu_1 \bar
\mu_2}{\sqrt{I(\bar \lambda_1, \bar \mu_1)} \sqrt{I(\bar
\lambda_2, \bar \mu_2)}}}= \varepsilon\, \ds{\frac{L \lambda_1
\lambda_2 + M (\lambda_1 \mu_2 + \mu_1 \lambda_2) + N \mu_1
\mu_2}{\sqrt{I(\lambda_1, \mu_1)} \sqrt{I(\lambda_2, \mu_2)}}} =
\varepsilon \, \zeta_{\,g_1,g_2}.$$ Consequently,
$\zeta_{\,g_1,g_2}$ is invariant (up to the orientation of the
tangent space or the normal space of the surface). \qed

\vskip 3mm
\begin{defn} \label{D:conjugate}
Two tangents $g_1: X_1 = \lambda_1 z_u + \mu_1 z_v$ and $g_2: X_2
= \lambda_2 z_u + \mu_2 z_v$ are said to be \emph{conjugate
tangents}   if $\zeta_{\,g_1,g_2} = 0$.
\end{defn}

Obviously, $\zeta_{\,g_1,g_2} = 0$ if and only if
$$L\lambda_1 \lambda_2 + M (\lambda_1 \mu_2 +\lambda_2 \mu_1) + N\mu_1 \mu_2 = 0.$$

The last formula gives us the idea to define second fundamental
form $II$ of the surface $M^2$ at $p \in M^2$ as follows. Let
$X=\lambda z_u+\mu z_v, \,\, (\lambda,\mu)\neq(0,0)$ be a tangent
vector at a point $p \in M^2$. Then
$$II(\lambda,\mu)=L\lambda^2+2M\lambda\mu+N\mu^2, \quad \lambda, \mu \in {\R}.$$
We have already seen in the proof of Proposition 2.1 that the
functions $L,M,N$ change in the same way as the coefficients
$E,F,G$  under any change of the parameters on $M^2$. If
$\{\widetilde{n}_1, \widetilde{n}_2\}$ is another normal frame
field of $M^2$, such that $\langle \widetilde{n}_1,
\widetilde{n}_1 \rangle = 1$, $\langle \widetilde{n}_1,
\widetilde{n}_2 \rangle = 0$, $\langle \widetilde{n}_2,
\widetilde{n}_2 \rangle = -1$, then
$$\begin{array}{l}
\vspace{2mm}
n_1 = \varepsilon'(\cosh \theta \, \widetilde{n}_1 + \sinh \theta \, \widetilde{n}_2);\\
\vspace{2mm} n_2 = \varepsilon'(\sinh \theta \, \widetilde{n}_1 +
\cosh \theta \, \widetilde{n}_2);
\end{array} \qquad \varepsilon'=\pm 1.$$ The relation between the corresponding functions
$c_{ij}^k$ and $\widetilde{c}_{ij}^k$, $i,j,k = 1,2$ is given by
the equalities
$$\begin{array}{l}
\vspace{2mm}
\widetilde{c}_{ij}^1 = \varepsilon'(\cosh \theta \, c_{ij}^1 - \sinh \theta \, c_{ij}^2);\\
\vspace{2mm} \widetilde{c}_{ij}^2 =\varepsilon'( - \sinh \theta \,
c_{ij}^1 + \cosh \theta \, c_{ij}^2);
\end{array}
\quad i,j = 1,2.$$ Thus,  $\widetilde{\Delta}_i =  \Delta_i$, $i =
1,2,3$, and $\widetilde{L} =  L, \; \widetilde{M} =  M, \;
\widetilde{N} =  N$. So, the functions $L$, $M$, $N$ do not depend
on the normal frame of the surface. Hence, the second fundamental
form $II$ is invariant up to the orientation of the tangent space
or the normal space of the surface.

As in the classical differential geometry of surfaces in $\R^3$
and in the same way as in the theory of surfaces on $\R^4$
\cite{GM3}, the second fundamental form $II$ determines conjugate
tangents at a point $p$ of $M^2$. The considerations above  show
that the conjugacy in terms of the second fundamental form is the
conjugacy defined by the invariant $\zeta_{\,g_1,g_2}$.

\vskip 2mm We shall assign two invariants $\nu_g$ and $\alpha_g$
to any tangent $g$ of the surface in the following way. Let $g: X
= \lambda z_u + \mu z_v$ be a tangent and $g^{\bot}$ be its
orthogonal tangent, determined by the vector
$$X^{\bot} = \ds{-\frac{F \lambda + G \mu}{W} \, z_u + \frac{E \lambda + F \mu}{W} \, z_v}.  \leqno{(2.8)}$$

We define
$$\nu_g = \zeta_{\,g,g}; \qquad \alpha_g = \zeta_{\,g, g^{\bot}}. $$
We call $\nu_g$ the \emph{normal curvature} of the tangent $g$,
and $\alpha_g$ - the \emph{geodesic torsion} of $g$.

Equalities (2.4) and (2.5) imply that
$$\nu_g = \ds{2\frac{S(\sigma_g(z_u), \sigma_g(z_v))}{W} = \frac{II(\lambda, \mu)}{I(\lambda,\mu)}}.$$
Hence, the normal curvature of the tangent $g$ is two times  the
oriented area of the parallelogram determined by the normal
vectors $\sigma_g(z_u)$ and $\sigma_g(z_v)$. The invariant $\nu_g$
is expressed by the first and the second fundamental forms of the
surface in the same way as the normal curvature of a tangent in
the theory of surfaces in $\R^3$.

Using (2.5) and (2.8) we get
$$\alpha_g = \ds{\frac{\lambda^2 (EM - FL) + \lambda \mu (EN - GL) + \mu^2(FN - GM)}{W I(\lambda, \mu)}}.$$
Hence, $\alpha_g$ is expressed by the coefficients of the first
and the second fundamental forms in the same way as the geodesic
torsion in the theory of surfaces in $\R^3$.

\vskip 2mm We define  asymptotic tangents and principal tangents
as follows:

\vskip 3mm
\begin{defn} \label{D:asymptotic tangent}
A tangent $g: X = \lambda z_u + \mu z_v$ is said to be
\emph{asymptotic} if it is self-conjugate, i.e. $\nu_g = 0$.
\end{defn}

\vskip 3mm
\begin{defn} \label{D:principal tangent}
A tangent $g: X = \lambda z_u + \mu z_v$ is said to be
\emph{principal} if it is perpendicular to its conjugate, i.e.
$\alpha_g = 0$.
\end{defn}

\vskip 3mm The equation for the asymptotic tangents at a point $p
\in M^2$ is
$$L\lambda^2 + 2M \lambda \mu + N\mu^2 = 0.$$

\vskip 2mm The equation for the principal tangents at a point $p
\in M^2$ is
$$\left|\begin{array}{cc}
E & F\\
[2mm] L & M \end{array}\right| \lambda^2+ \left|\begin{array}{cc}
E & G\\
[2mm] L & N \end{array}\right| \lambda \mu+
\left|\begin{array}{cc}
F & G\\
[2mm] M & N \end{array}\right| \mu^2=0.$$

A line $c: u=u(q), \; v=v(q); \; q\in J \subset \R$ on $M^2$ is
said to be an \emph{asymptotic line}, respectively a
\textit{principal line}, if its tangent at any point is
asymptotic, respectively  principal. The surface $M^2$ is
parameterized by principal lines if and only if $F=0, \,\, M=0.$

\section{Weingarten map of a spacelike surface in $\R^4_1$} \label{S:Weingarten}

The second fundamental form $II$ determines a map of
Weingarten-type $\gamma: T_pM^2 \rightarrow T_pM^2$ at any point
of $M^2$ in the standard way:
$$\begin{array}{l}
\vspace{2mm}
\gamma(z_u)=\gamma_1^1z_u+\gamma_1^2z_v,\\
\vspace{2mm} \gamma(z_v)=\gamma_2^1z_u+\gamma_2^2z_v,
\end{array}$$
where
$$\displaystyle{\gamma_1^1=\frac{FM-GL}{EG-F^2}, \quad
\gamma_1^2 =\frac{FL-EM}{EG-F^2}}, \quad
\displaystyle{\gamma_2^1=\frac{FN-GM}{EG-F^2}, \quad
\gamma_2^2=\frac{FM-EN}{EG-F^2}}.$$ The linear map $\gamma$ is
invariant under changes of the parameters of the surface and
changes of the normal frame field. Hence the following statement
holds.

\begin{lem}\label{L:Lemma-invariant}
The functions
$$k := \det \gamma = \frac{LN - M^2}{EG - F^2}, \qquad
\varkappa :=-\frac{1}{2}\,{\rm tr}\, \gamma
=\frac{EN+GL-2FM}{2(EG-F^2)}$$ are invariants of the surface
$M^2$.
\end{lem}

Now we shall prove

\vskip 2mm
\begin{prop}
The function $\varkappa$ is the curvature of the normal connection
of $M^2$.
\end{prop}

\noindent \emph{Proof:} Let $D$ be the normal connection of $M^2$.
For any tangent vector fields $x,y$ and any normal vector field
$n$ we have the standard decomposition
$$\nabla'_xn = -A_n(x) + D_xn,$$
where $\langle A_n(x), y\rangle = \langle \sigma(x,y),n\rangle$.

The curvature tensor $R^{\bot}$ of the normal connection $D$ is
given by $R^{\bot}(x,y)n = D_xD_yn - D_yD_xn - D_{[x,y]}n$. Then
the curvature of the normal connection at a point $p \in M^2$ is
defined by $\langle R^{\bot}(x,y)n_2, n_1\rangle,$ where
$\{x,y,n_1,n_2\}$ is a right oriented orthonormal quadruple.

Without loss of generality we assume that $F = 0$ and denote the
unit vector fields $\displaystyle{x=\frac{z_u}{\sqrt E}, \;
y=\frac{z_v}{\sqrt G}}$. Then we have
$$\begin{array}{l}
\vspace{2mm}
\sigma(x,x)=\displaystyle{\frac{c_{11}^1}{E}\;\;n_1 \;-\;\frac{c_{11}^2}{E}\;\;n_2},\\
\vspace{2mm}
\sigma(x,y)= \displaystyle{\frac{c_{12}^1}{\sqrt{EG}}\,n_1 -\frac{c_{12}^2}{\sqrt{EG}}\,n_2,}\\
\vspace{2mm} \sigma(y,y) =\displaystyle{\frac{c_{22}^1}{G}\;\;n_1
\;-\;\;\frac{c_{22}^2}{G}\;\;n_2}.
\end{array}$$
Hence,
$$\begin{array}{ll}
\vspace{2mm} A_1(x)=\displaystyle{\frac{c_{11}^1}{E}\;\;x
\;+\;\frac{c_{12}^1}{\sqrt{EG}}\;\;y},
\qquad & A_2(x)=\displaystyle{\frac{c_{11}^2}{E}\;\;x \;+\;\frac{c_{12}^2}{\sqrt{EG}}\;\;y},\\
\vspace{2mm} A_1(y)= \displaystyle{\frac{c_{12}^1}{\sqrt{EG}}\;\;x
+ \frac{c_{22}^1}{G}\;\;y,} \qquad & A_2(y)=
\displaystyle{\frac{c_{12}^2}{\sqrt{EG}}\;\;x +
\frac{c_{22}^2}{G}\;\;y.}
\end{array} \leqno{(3.1)}$$
Using (3.1) we calculate
$$\begin{array}{ll}
\vspace{2mm} (A_2 \circ A_1 - A_1 \circ A_2) (x) =
\ds{\left(\frac{c_{11}^1 c_{12}^2 - c_{11}^2 c_{12}^1}{E
\sqrt{EG}} + \frac{c_{12}^1 c_{22}^2 - c_{12}^2 c_{22}^1}{G
\sqrt{EG}}\right) y} =
\ds{\frac{EN + GL}{2EG}\, y;}\\
\vspace{2mm} (A_2 \circ A_1 - A_1 \circ A_2) (y) = -
\ds{\left(\frac{c_{11}^1 c_{12}^2 - c_{11}^2 c_{12}^1}{E
\sqrt{EG}} + \frac{c_{12}^1 c_{22}^2 - c_{12}^2 c_{22}^1}{G
\sqrt{EG}}\right) x} = - \ds{\frac{EN + GL}{2EG}\, x.}
\end{array}$$
Hence,
$$\begin{array}{ll}
\vspace{2mm}
(A_2 \circ A_1 - A_1 \circ A_2) (x)= \varkappa \,y;\\
\vspace{2mm} (A_2 \circ A_1 - A_1 \circ A_2) (y)= - \varkappa \,x.
\end{array} \leqno{(3.2)}$$

Note that $A_2 \circ A_1 - A_1 \circ A_2$ is an invariant
skew-symmetric operator in the tangent space, i.e. it does not
depend on the choice of the orthonormal tangent frame field
$\{x,y\}$.

Since the curvature tensor $R'$ of the connection $\nabla'$ is
zero, we have
$$\nabla'_x \nabla'_y n_1 - \nabla'_y \nabla'_x n_1 - \nabla'_{[x,y]} n_1 = 0.$$
Therefore the tangent component and the normal component of
$R'(x,y)n_1$ are both zero. The normal component is $D_xD_yn_1 -
D_yD_xn_1 - D_{[x,y]}n_1 - \sigma(x,A_1y) + \sigma(y,A_1x).$
Hence,
$$D_xD_yn_1 - D_yD_xn_1 - D_{[x,y]}n_1 = \sigma(x,A_1y) - \sigma(y,A_1x). \leqno{(3.3)}$$
The left-hand side of (3.3) is $R^{\bot}(x,y)n_1$. Then
$$\langle R^{\bot}(x,y)n_1, n_2 \rangle = \langle \sigma(x,A_1y), n_2 \rangle
- \langle \sigma(y,A_1x), n_2 \rangle = \langle (A_2 \circ A_1 -
A_1 \circ A_2) (y), x \rangle.$$ Using (3.2) we obtain that
$\langle R^{\bot}(x,y)n_1, n_2 \rangle = - \varkappa$. Since
$\langle R^{\bot}(x,y)n_1, n_2 \rangle = - \langle
R^{\bot}(x,y)n_2, n_1 \rangle$, we get
$$\langle R^{\bot}(x,y)n_2, n_1 \rangle = \varkappa.$$
The last equality implies that $\varkappa$ is the curvature of the
normal connection. \qed

\vskip 3mm

The characteristic equation of the Weingarten map $\gamma$ is
$$\nu^2 + 2 \varkappa \, \nu + k = 0.$$
Since $\gamma$ is a symmetric linear operator, the following
inequality holds:
$$\varkappa^2 - k \geq 0.$$
Moreover, the equality $\varkappa^2 - k  = 0$ is equivalent to the
conditions
$$L = \rho E, \quad M = \rho F, \quad N = \rho G, \qquad \rho\in {\R}.$$

Obviously, the following equivalence at a point $p \in M^2$ holds:
$$L = M = N = 0 \quad \iff \quad k = \varkappa = 0.$$

As in the theory of surfaces in $\R^3$ and $\R^4$, the invariants
$k$ and $\varkappa$ divide the points of $M^2$ into four types. A
point $p \in M^2$ is said to be: \vskip 1mm \emph{flat}, \; if \;
$k = \varkappa = 0$; \vskip 1mm \emph{elliptic}, \; if \; $k > 0$;
\vskip 1mm \emph{parabolic}, \; if \; $k = 0$, \; $\varkappa \neq
0$; \vskip 1mm \emph{hyperbolic}, \; if \; $k < 0$.

\vskip 1mm Spacelike surfaces consisting of flat points will be
considered in Section \ref{S:flat points}. Further in this section
we shall consider surfaces free of flat points, i.e. $(L,M,N) \neq
(0,0,0)$.

\vskip 2mm We note that a spacelike surface $M^2$ has two families
of orthogonal asymptotic lines if and only if $M^2$ is  of flat
normal connection.

\vskip 2mm Let $\ds{H=\frac{1}{2}\, (\sigma(x, x)+\sigma(y,y))}$
be the normal mean curvature vector field. We recall that a
surface $M^2$ is said to be \textit{minimal} if the mean curvature
vector $H = 0$. The minimal surfaces are characterized in terms of
the invariants $k$ and $\varkappa$ by the following

\vskip 2mm
\begin{prop}\label{P:minimal}
Let $M^2$ be a spacelike surface in $\R^4_1$ free of flat points.
Then $M^2$ is minimal if and only if
$$\varkappa^2 - k = 0.$$
\end{prop}

\vskip 2mm \noindent \emph{Proof:} Without loss of generality we
assume that $F = 0$ and denote the unit vector fields
$\displaystyle{x=\frac{z_u}{\sqrt E}, \; y=\frac{z_v}{\sqrt G}}$.
Then we have
$$\begin{array}{l}
\vspace{2mm} \nabla'_xx=\quad \quad \quad
\gamma_1\,y+\;\;\displaystyle{\frac{c_{11}^1}{E}\;\;n_1
\;-\;\frac{c_{11}^2}{E}\;\;n_2},\\
\vspace{2mm} \nabla'_xy=-\gamma_1\,x\quad \quad \;
+\displaystyle{\frac{c_{12}^1}{\sqrt{EG}}\,n_1
-\frac{c_{12}^2}{\sqrt{EG}}\,n_2,}\\
\vspace{2mm} \nabla'_yx=\quad\quad \; \; -\gamma_2\,y
+\displaystyle{\frac{c_{12}^1}{\sqrt{EG}}\,n_1
-\frac{c_{12}^2}{\sqrt{EG}}\,n_2,}\\
\vspace{2mm} \nabla'_yy=\;\;\gamma_2\,x
\quad\quad\quad+\;\;\displaystyle{\frac{c_{22}^1}{G}\;\;n_1
\;-\;\;\frac{c_{22}^2}{G}\;\;n_2}.\\
\end{array}$$

I. Let $\ds{H=\frac{1}{2}\, (\sigma(x, x)+\sigma(y,y))=0}$. Then
$c_{22}^1 = \ds{-\frac{G}{E} \,c_{11}^1}$, $c_{22}^2 =
\ds{-\frac{G}{E} \,c_{11}^2}$,
 and hence
$$\Delta_2 = \left|%
\begin{array}{cc}
\vspace{2mm}
  c_{11}^1 & c_{22}^1 \\
  c_{11}^2 & c_{22}^2 \\
\end{array}%
\right|=0, \quad \quad \frac{\Delta_3}{G}=\frac{\Delta_1}{E}.$$
Therefore
$$L = \rho E, \quad M = \rho F, \quad N = \rho G,$$
where $\rho$ is a function on $M^2$. Hence $\varkappa^2 - k = 0$.

II. Let $\varkappa^2 - k = 0$. Then
$$L = \rho E, \quad M = \rho F, \quad N = \rho G; \quad \rho \neq 0.$$
The condition $F=0$ implies that $M=0$. Then
$\displaystyle{\left|%
\begin{array}{cc}
\vspace{2mm}
  c_{11}^1 & c_{22}^1 \\
  c_{11}^2 & c_{22}^2 \\
\end{array}%
\right|=0}$ and $c_{22}^1=\widetilde{\rho} c_{11}^1, \;
c_{22}^2=\widetilde{\rho} c_{11}^2$. Further, the equality
$\displaystyle{\frac{L}{E}=\frac{N}{G}}$ implies that
$\displaystyle{\widetilde{\rho}=-\frac{G}{E}}$. Hence $\tr \,
\sigma = 0$, i.e. $H=0$. \qed \vskip 2mm

Let us note that the spacelike surfaces consisting of
''umbilical'' points in $\R^4_1$ are exactly the minimal surfaces.

\vskip 3mm We shall  characterize the minimal surfaces and the
surfaces with flat normal connection in terms of  a geometric
figure in the tangent space at any point of a spacelike surface.

The normal curvatures  of the principal tangents are said to be
\textit{principal normal curvatures} of $M^2$. If a point $p \in
M^2$ is ''non-umbilical'', i.e. $\varkappa^2 - k>0$, we can assume
that $(u,v)$ are principal parameters ($F=M=0$). The  principal
normal curvatures are $\nu' = \ds{\frac{L}{E}}$, $\nu'' =
\ds{\frac{N}{G}}$  and the invariants $k$ and $\varkappa$ of $M^2$
are expressed by the principal normal curvatures  as follows:
$$k = \nu' \nu''; \qquad \varkappa = \frac{\nu' + \nu''}{2}. \leqno{(3.4)}$$
If $p \in M^2$ is an ''umbilical'' point, i.e. $\varkappa^2 -
k=0$, then all tangents at $p$ are principal with one and the same
normal curvature $\nu$. Then formulas (3.4) are also valid under
the assumption $\nu'=\nu''=\nu$.

\vskip 1mm

Similarly to the theory of surfaces in $\R^3$ and $\R^4$, we
consider the indicatrix $\chi$ in the tangent space $T_pM^2$ at an
arbitrary point $p$ of $M^2$, defined by
$$\chi : \nu' X^2 + \nu'' Y^2 = \varepsilon, \qquad \varepsilon = \pm 1.$$

If $p$ is an elliptic point ($k > 0$), then the indicatrix $\chi$
is an ellipse. The axes of $\chi$ are collinear with the principal
tangents at the point $p$, and the lengths of the axes are
$\ds{\frac{2}{\sqrt{|\nu'|}}}$ and
$\ds{\frac{2}{\sqrt{|\nu''|}}}$.

If $p$ is a  hyperbolic point ($k < 0$), then the indicatrix
$\chi$ consists of two hyperbolas. For the sake of simplicity we
say that $\chi$ is a hyperbola. The axes of $\chi$ are collinear
with the principal tangents, and the lengths of the axes are
$\ds{\frac{2}{\sqrt{|\nu'|}}}$ and
$\ds{\frac{2}{\sqrt{|\nu''|}}}$.

If $p$ is a parabolic point ($k = 0$), then the indicatrix $\chi$
consists of two straight lines parallel to the principal tangent
with non-zero normal curvature.

The following statement holds:

\begin{prop}\label{P:conjugacy}
Two tangents $g_1$ and $g_2$ are conjugate tangents of $M^2$ if
and only if $g_1$ and $g_2$ are conjugate with respect to the
indicatrix $\chi$.
\end{prop}

\vskip 2mm Now we shall characterize the minimal surfaces and the
surfaces with flat normal connection in terms of the tangent
indicatrix of the surface.

\begin{prop}\label{P:Minimal-circle}
Let $M^2$ be a spacelike surface in $\R^4_1$ free of flat points.
Then $M^2$ is minimal if and only if at each point of $M^2$ the
tangent indicatrix $\chi$ is a circle.
\end{prop}

\noindent {\it Proof:} Let $M^2$ be a spacelike surface in
$\R^4_1$ free of flat points. From equalities (3.4) it follows
that
$$\varkappa^2 - k = \ds{\left(\frac{\nu' - \nu''}{2}\right)^2}.$$
Obviously $\varkappa^2 - k = 0$ if and only if $\nu' = \nu''$.
Applying Proposition \ref{P:minimal}, we get that $M^2$ is minimal
if and only if  $\chi$ is a circle. \qed

\begin{prop}\label{P:Flat normal-Lorentz circle}
Let $M^2$ be a spacelike surface in $\R^4_1$ free of flat points.
Then $M^2$ is a surface with flat normal connection if and only if
at each point of $M^2$ the tangent indicatrix $\chi$ is a
rectangular hyperbola (a Lorentz circle).
\end{prop}

\noindent {\it Proof:} Let $M^2$ be a spacelike surface in
$\R^4_1$ free of flat points. From (3.4) it follows that
$\varkappa = 0$ if and only if $\nu'' = -\nu'$.

If $M^2$ is a surface with flat normal connection, then $k < 0$,
and hence $\chi$ is a hyperbola. From $\nu'' = -\nu'$ it follows
that the semi-axes of $\chi$ are equal to
$\ds{\frac{1}{\sqrt{|\nu'|}}}$, i.e. $\chi$ is a rectangular
hyperbola.

Conversely, if $\chi$ is a rectangular hyperbola, then $\nu'' =
-\nu'$, which implies that $M^2$ is a surface with flat normal
connection. \qed

\vskip 3mm The minimal surfaces and the surfaces with flat normal
connection can also be characterized in terms of the ellipse of
normal curvature.

The notion of the ellipse of normal curvature of a surface in
space forms was introduced by Moore and Wilson \cite{MW1, MW2}.
The ellipse of normal curvature associated to the second
fundamental form of a spacelike surface in $\R^4_1$ was first
considered in \cite{IPF}. The \textit{ellipse of normal curvature}
at a point $p$ of a surface $M^2$  is the ellipse in the normal
space at the point $p$ given by $\{\sigma(x,x): \, x \in T_pM^2,
\, \langle x,x \rangle = 1\}$.
 Let $\{x,y\}$ be an orthonormal base of the
tangent space $T_pM^2$ at $p$. Then, for any $v = \cos \psi \, x +
\sin \psi \, y$, we have
$$\sigma(v, v) = H + \ds{\cos 2\psi  \, \frac{\sigma(x,x) - \sigma(y,y)}{2}
+ \sin 2 \psi  \, \sigma(x,y)}, \leqno{(3.5)}$$ where $H $  is the
mean curvature vector of $M^2$ at $p$. So, when $v$ goes once
around the unit tangent circle, the vector $\sigma(v,v)$ goes
twice around the ellipse centered at $H$.
 The vectors $\ds{\frac{\sigma(x,x) - \sigma(y,y)}{2}}$ \, and
$\sigma(x,y)$ determine conjugate directions of the ellipse.

\vskip 2mm Obviously, $M^2$ is minimal if and only if for each
point $p \in M^2$ the ellipse of curvature is centered at $p$. We
shall give a characterization of the surfaces with flat normal
connection in terms of the ellipse of normal curvature.

\begin{lem}\label{L:Flat normal-ellipse}
Let $M^2$ be a spacelike surface in $\R^4_1$ free of flat points
and $x$, $y$ be principal tangents. Then $M^2$ is a surface with
flat normal connection if and only if $\sigma(x,x) = \sigma(y,y)$.
\end{lem}

\vskip 2mm \noindent {\it Proof:} Let $M^2$ be a surface in
$\R^4_1$ free of flat points, and parameterized by principal
parameters, i.e. $F = M =0$. Then $\varkappa = \ds{\frac{EN +
GL}{2EG}}$. Denote $\displaystyle{x=\frac{z_u}{\sqrt E}, \;
y=\frac{z_v}{\sqrt G}}$. Since $M=0$, we have
$\displaystyle{\left|%
\begin{array}{cc}
\vspace{2mm}
  c_{11}^1 & c_{22}^1 \\
  c_{11}^2 & c_{22}^2 \\
\end{array}%
\right|=0}$ and hence $c_{22}^1=\rho \,c_{11}^1, \;
c_{22}^2=\rho\, c_{11}^2$. Then $\Delta_3 = - \rho\,\Delta_1$ and
$N= - \rho \, L$.

I. Let $M^2$ be of flat normal connection, i.e. $\varkappa = 0$.
Then, $N = \ds{-\frac{G}{E} L}$, and hence, $\rho =
\ds{\frac{G}{E}}$, which implies that
 $c_{22}^1=\ds{\frac{G}{E}} \,c_{11}^1, \; c_{22}^2=\ds{\frac{G}{E}}\, c_{11}^2$.
Consequently,
$$\sigma(y,y) = \displaystyle{\frac{c_{22}^1}{G}\,n_1
-\frac{c_{22}^2}{G}\,n_2} = \displaystyle{\frac{c_{11}^1}{E}\,n_1
-\frac{c_{11}^2}{E}\,n_2} = \sigma(x,x).$$

II. Let $\sigma(x,x) = \sigma(y,y)$. Then
$\ds{\frac{c_{22}^1}{G}=\frac{c_{11}^1}{E}}; \;
\ds{\frac{c_{22}^2}{G}= \frac{c_{11}^2}{E}}$. Using that
$c_{22}^1=\rho \,c_{11}^1, \; c_{22}^2=\rho\, c_{11}^2$, we get
$\rho = \ds{\frac{G}{E}}$, and hence $N = -\ds{\frac{G}{E} L}$.
The last equality implies $\varkappa = 0$, i.e. $M^2$ is a surface
of flat normal connection. \qed

\vskip 3mm
\begin{prop}\label{P:Flat normal-ellipse}
Let $M^2$ be a spacelike surface in $\R^4_1$ free of flat points.
Then $M^2$ is a surface of flat normal connection if and only if
for each point $p \in M^2$ the ellipse of normal curvature is a
line segment, which is not collinear with the mean curvature
vector field.
\end{prop}

\vskip 2mm \noindent {\it Proof:} I. Let $M^2$ be a surface of
flat normal connection. According to Lemma \ref{L:Flat
normal-ellipse} we have $\sigma(x,x) - \sigma(y,y) = 0$. Then for
any $v = \cos \psi \, x + \sin \psi \, y$, we get $\sigma(v, v) =
H + \sin 2 \psi  \, \sigma(x,y)$. Hence, when $v$ goes once around
the unit tangent circle, the vector $\sigma(v,v)$ goes twice along
the line segment collinear with $\sigma(x,y)$ and centered at $H$.
If we assume that $\sigma(x,y)$ is collinear with $H =
\sigma(x,x)$, we get $L = M = N =0$, which contradicts the
condition that $M^2$ is free of flat points.

II. Let $M^2$ be a surface in $\R^4_1$ free of flat points such
that for each point $p \in M^2$ the ellipse of normal curvature is
a line segment, which is not collinear with $H$. Without loss of
generality we assume that $M^2$ is parameterized by principal
parameters, and hence $c_{22}^1=\rho \,c_{11}^1, \;
c_{22}^2=\rho\, c_{11}^2$. Then
$$\begin{array}{l}
\vspace{2mm}
\sigma(x,x) - \sigma(y,y) = (1 - \rho) (c_{11}^1\,n_1 + c_{11}^2\,n_2);\\
\vspace{2mm} \sigma(x,y) = c_{12}^1\,n_1 + c_{12}^2\,n_2.
\end{array} \leqno{(3.6)}$$
Since the ellipse of normal curvature is a line segment, having in
mind (3.5), we get one of the following  possibilities:

(a) $\sigma(x,x) - \sigma(y,y)$ is collinear with  $\sigma(x,y)$.
In this case from (3.6) we get $ c_{12}^1 = \widetilde{\rho}
c_{11}^1; \,c_{12}^2 = \widetilde{\rho}  c_{11}^2$, which implies
$L=M=N=0$, a contradiction.

(b) $\sigma(x,y) =0$, which implies $c_{12}^1 = c_{12}^2 = 0$, and
hence $L=M=N=0$, a contradiction.

(c) $\sigma(x,x) - \sigma(y,y) = 0$. i.e. $\sigma(x,x) =
\sigma(y,y)$. Applying Lemma \ref{L:Flat normal-ellipse}, we get
that $M^2$ is a surface of flat normal connection. \qed

\section{Spacelike surfaces consisting of flat points}\label{S:flat points}

In this section we consider spacelike surfaces consisting of flat
points, i.e. surfaces satisfying the conditions
$$k(u,v)=0, \quad \varkappa(u,v)=0, \qquad (u,v) \in \mathcal D,$$
or equivalently $L(u,v)=0, \,\, M(u,v)=0, \,\,N(u,v) = 0, \,\,
(u,v) \in \mathcal D.$

We shall give a local geometric description of those such surfaces
whose mean curvature vector $H$ at any point is:

(1) a non-zero spacelike vector, i.e. $\langle H,H \rangle >0$, or

(2) a timelike vector, i.e. $\langle H,H \rangle <0$.

\begin{thm}\label{T:Main Theorem}
Let  $M^2$  be a spacelike surface in $\R^4_1$ consisting of flat
points and the mean curvature vector at any point is a non-zero
spacelike vector or timelike vector. Then either $M^2$ lies in a
hyperplane of $\R^4_1$  or $M^2$ is part of a developable ruled
surface in $\R^4_1$.

\end{thm}

\vskip 2mm \noindent \emph{Proof:} Let $M^2: z = z(u,v), \, \,
(u,v) \in {\mathcal D}$ (${\mathcal D} \subset \R^2$) be a
spacelike surface in $\R^4_1$ whose mean curvature vector at any
point is a non-zero spacelike vector or  timelike vector and
$L(u,v)= M(u,v)= N(u,v) = 0, \,\, (u,v) \in \mathcal D$. For the
sake of simplicity, we assume that the parametrization of $M^2$ is
orthogonal, i.e. $F=0$. Denote the unit vector fields
$\displaystyle{x=\frac{z_u}{\sqrt E}, \; y=\frac{z_v}{\sqrt G}}$.
The conditions $L = M = N = 0$ imply that
$${\rm rank}\left(
\begin{array}{ccc}
c_{11}^1 & c_{12}^1 & c_{22}^1\\
[2mm] c_{11}^2 & c_{12}^2 & c_{22}^2
\end{array}\right)=1 $$
and the vectors $\sigma (x,x), \; \sigma(x,y), \; \sigma(y,y)$ are
collinear. Let $n$ be a unit normal vector field of $M^2$, which
is collinear with $\sigma (x,x), \; \sigma(x,y)$, and
$\sigma(y,y)$. Hence, $n$ is collinear with the mean curvature
vector field $H$. We have the following possibilities:

(1) $n$ is spacelike, i.e. $\langle n,n \rangle =1$.

(2) $n$ is timelike, i.e. $\langle n,n \rangle = -1$.

First we shall consider the case $\langle n, n\rangle=1$. Denote
by $l$ the unit normal vector field such that $\{x,y,n,l\}$ is a
positively oriented orthonormal frame field in $\R^4_1$ (hence
$\langle l, l\rangle=-1$). It is clear that the normal vector
fields $n,l$ are determined  up to a sign. Then we have the
following derivative formulas of $M^2$:
$$\begin{array}{ll}
\vspace{2mm} \nabla'_xx=\quad \quad \quad \gamma_1\,y+\;\nu_1\,n,&
\qquad \quad
\nabla'_xn=-\nu_1\,x-\lambda\,y \quad \, \quad- \beta_1\,l,\\
\vspace{2mm} \nabla'_xy=-\gamma_1\,x\quad \quad \; + \;
\lambda\,n, & \qquad \quad
\nabla'_yn=-\lambda \, x -\nu_2\,y\quad \quad \,- \beta_2 \,l,\\
\vspace{2mm} \nabla'_yx=\quad\quad \; -\gamma_2\,y  + \;\lambda\,
n, & \qquad \quad
\nabla'_xl=\quad\quad \quad \quad \quad-\beta_1\,n,\\
\vspace{2mm} \nabla'_yy=\;\;\gamma_2\,x \quad\quad
\;\;+\;\nu_2\,n, & \qquad \quad \nabla'_yl=\quad\quad \quad \quad
\quad-\beta_2\,n,
\end{array}\leqno(4.1)$$
where $\nu_1, \, \nu_2,\, \lambda,\,\beta_1,\, \beta_2,
\gamma_1,\, \gamma_2$ are functions on $M^2$.

The mean curvature vector field is $H = \ds{\frac{\nu_1 +
\nu_2}{2} \,n}$. The Gauss curvature $K$ of $M^2$ is expressed by
$$K = \nu_1\,\nu_2-\lambda^2.\leqno(4.2)$$

Since the curvature tensor $R'$ of the connection $\nabla'$ is
zero, then the equality  $R'(x,y,l) = 0$ together with  (4.1)
imply that
$$\begin{array}{ll}
\vspace{2mm}
\nu_1 \,\beta_2 - \lambda \,\beta_1 =0;\\
\vspace{2mm} -\lambda \,\beta_2 + \nu_2 \,\beta_1 =0.
\end{array}\leqno(4.3)$$
We have two subcases:

\hskip 10mm (a) $\beta_1 = \beta_2 = 0$ for all $(u,v) \in
\mathcal D$. Then from equalities (4.1) it follows that
$\nabla'_xl=  0;\,\, \nabla'_yl=0$, and hence $l={\rm const}$.
Consequently, $M^2$ lies in a hyperplane ${\E}^3$ of $\R^4_1$
orthogonal to $l$.

\vskip 1mm \hskip 10mm (b) There exists a point $(u_0,v_0) \in
\mathcal D$ such that $\beta_1^2(u_0,v_0) +  \beta_2^2(u_0,v_0)
\neq 0$. Hence, there exists a neighborhood ${\mathcal D}_0
\subset \mathcal D$ such that $({\beta_1^2 +
\beta_2^2})_{|{\mathcal D}_0} \neq 0$. Then, equalities (4.3)
imply that $\nu_1 \, \nu_2 - \lambda^2 = 0$ for $(u,v) \in
{\mathcal D}_0$. If $\nu_1 = \nu_2 = 0$, then  $H=0$, which
contradicts the assumption that $\langle H,H \rangle >0$. So we
assume that there exists a neighborhood $\widetilde{\mathcal D}
\subset {\mathcal D}_0$ such that ${\nu_2}_{|\widetilde{\mathcal
D}} \neq 0$ (or ${\nu_1}_{|\widetilde{\mathcal D}} \neq 0$). We
consider the surface $\widetilde{M}^2 = M^2_{|\widetilde{\mathcal
D}}$, which is a surface with zero Gauss curvature in view of
(4.2).

Let $\{\overline{x},\overline{y}\}$ be the orthonormal tangent
frame field of $\widetilde{M}^2$, defined by
$$\begin{array}{l}
\vspace{2mm}
\overline{x} = \cos \varphi \, x + \sin \varphi \, y;\\
\vspace{2mm} \overline{y} = - \sin \varphi \, x +  \cos \varphi
\,y,
\end{array}$$
where $\tan \varphi = \displaystyle{- \frac{\lambda}{\nu_2}}$.
Then $\sigma(\overline{x}, \overline{x}) = 0$, \,
$\sigma(\overline{x}, \overline{y}) = 0$. So, formulas (4.1) take
the form
$$\begin{array}{ll}
\vspace{2mm} \nabla'_{\overline{x}}\,\overline{x} = \quad \quad
\quad\;\; \overline{\gamma}_1\,\overline{y},&
\qquad \quad \nabla'_{\overline{x}}\,n = \quad \quad \quad \quad \quad \quad \,\,-\overline{\beta}_1\,l,\\
\vspace{2mm} \nabla'_{\overline{x}}\,\overline{y} =
-\overline{\gamma}_1\,\overline{x},&\qquad \quad
\nabla'_{\overline{y}}\,n = - \overline{\nu}_2\, \overline{y} \quad\quad \quad \;\;-\overline{\beta}_2\,l,\\
\vspace{2mm} \nabla'_{\overline{y}}\,\overline{x} = \quad \quad
\quad - \overline{\gamma}_2\,\overline{y},&
\qquad \quad\nabla'_{\overline{x}}\,l = \quad\quad \quad -\overline{\beta}_1\,n,\\
\vspace{2mm} \nabla'_{\overline{y}}\,\overline{y} =
\overline{\gamma}_2\,\overline{x} \quad \quad \quad \quad +
\overline{\nu}_2\, n,& \qquad \quad \nabla'_{\overline{y}}\,l =
\quad\quad \quad -\overline{\beta}_2\,n,
\end{array}$$
where $\overline{\nu}_2 \neq 0$.

Now the equalities $R'(\overline{x},\overline{y},n) = 0$ and
$R'(\overline{x},\overline{y},l) = 0$ imply that
$$\overline{\gamma}_1 = 0, \quad \quad \overline{\beta}_1 = 0.$$
Hence,
$$\begin{array}{ll}
\vspace{2mm}
\nabla'_{\overline{x}}\,\overline{x} = 0, & \quad \quad \nabla'_{\overline{x}}\,n = 0,\\
\vspace{2mm} \nabla'_{\overline{x}}\,\overline{y} = 0, & \quad
\quad \nabla'_{\overline{x}}\,l = 0.
\end{array}$$

Let $p= z(\overline{u}_0,\overline{v}_0), \,\,
(\overline{u}_0,\overline{v}_0) \in \widetilde{\mathcal D}$ be an
arbitrary point of  $\widetilde{M}^2$ and $c_1: z(\overline{u}) =
z(\overline{u},\overline{v}_0)$ be the integral curve of the
vector field $\overline{x}$, passing through $p$. It follows from
$\nabla'_{\overline{x}}\,\overline{x} = 0$  that $c_1$ is
contained in a straight line. Hence, $\widetilde{M}^2$ lies on a
one-parameter family of straight lines, i.e. $\widetilde{M}^2$
lies on a ruled surface. Moreover, since
$\nabla'_{\overline{x}}\,n = 0$ and $\nabla'_{\overline{x}}\,l =
0$ then the normal space $\span \{n, l\}$ of $\widetilde{M}^2$ is
constant at the points of $c_1$ and hence, the tangent space
$\span\{\overline{x},\overline{y}\}$ of $\widetilde{M}^2$ at the
points of $c_1$ is one and the same. Consequently,
$\widetilde{M}^2$ is part of a developable surface.

\vskip 2mm Now we shall consider the case $\langle n, n\rangle =
-1$. Denote by $b$ the unit normal vector field such that
$\{x,y,b,n\}$ is a positively oriented orthonormal frame field in
$\R^4_1$ (hence $\langle b, b\rangle=1$). The normal vector fields
$b,n$ are determined  up to a sign. In this case we have the
following derivative formulas of $M^2$:
$$\begin{array}{ll}
\vspace{2mm} \nabla'_xx=\quad \quad \quad \gamma_1\,y- \nu_1\,n,&
\qquad \quad
\nabla'_xb=  \qquad \qquad \qquad \qquad     - \beta_1\,n,\\
\vspace{2mm} \nabla'_xy=-\gamma_1\,x\quad \quad \; - \lambda\,n, &
\qquad \quad
\nabla'_yb= \qquad \qquad\qquad \qquad - \beta_2 \,n,\\
\vspace{2mm} \nabla'_yx=\quad\quad\; -\gamma_2\,y  - \lambda\, n,
& \qquad \quad
\nabla'_xn=-\nu_1\,x-\lambda\,y -\beta_1\,b,\\
\vspace{2mm} \nabla'_yy=\;\;\gamma_2\,x \quad\quad \;\;-\nu_2\,n,
& \qquad \quad \nabla'_yn=-\lambda \, x -\nu_2\,y -\beta_2\,b,
\end{array}\leqno(4.4)$$
where $\nu_1, \, \nu_2,\, \lambda,\,\beta_1,\, \beta_2,
\gamma_1,\, \gamma_2$ are functions on $M^2$.

The mean curvature vector field is $H = \ds{- \, \frac{\nu_1 +
\nu_2}{2} \;n}$. The Gauss curvature $K$ of $M^2$ is $K = -
\,\nu_1\,\nu_2+\lambda^2$. As in the previous case, using that
$R'(x,y,b) = 0$ from equalities (4.4) we obtain equalities (4.3)
which imply that there are two subcases:

\hskip 10mm (a) $\beta_1 = \beta_2 = 0$ for all $(u,v) \in
\mathcal D$. Then $\nabla'_xb=  0;\,\, \nabla'_yb=0$, and hence
$b={\rm const}$. Consequently, $M^2$ lies in a hyperplane ${\E}^3$
of $\R^4_1$ orthogonal to $b$.

\vskip 1mm \hskip 10mm (b) There exists a point $(u_0,v_0) \in
\mathcal D$ such that $\beta_1^2(u_0,v_0) +  \beta_2^2(u_0,v_0)
\neq 0$. In the same way as in the first case we obtain  that in a
neighborhood $\widetilde{\mathcal D} \subset \mathcal D$ the
surface $\widetilde{M}^2 = M^2_{|\widetilde{\mathcal D}}$ is part
of a developable surface. \qed

\section{Spacelike surfaces whose mean curvature vector at any point is a
non-zero spacelike vector}\label{S:spacelikeH}

Let $M^2$ be a spacelike surface parameterized by principal lines
and $\displaystyle{x=\frac{z_u}{\sqrt E}, \; y=\frac{z_v}{\sqrt
G}}$. The equality $M = 0$ implies that the normal vector fields
$\sigma(x,x)$ and $\sigma(y,y)$ are collinear. Hence, there exists
a geometrically determined normal frame field $n$, such that
$\sigma(x,x)$ and $\sigma(y,y)$ are collinear with $n$. Then we
have the following formulas:
$$\begin{array}{l}
\vspace{2mm}
\sigma(x,x) = \nu_1\,n; \\
\vspace{2mm} \sigma(y,y) = \nu_2\,n,
\end{array}$$
where $\nu_1, \nu_2$ are invariant functions. The mean curvature
vector field is expressed as follows:
$$H = \ds{\frac{\nu_1 + \nu_2}{2} \,n}.$$

Let $M^2$ be free of minimal points, i.e. $H \neq 0$ at each point
of $M^2$. We have the following possibilities for the mean
curvature vector field:
\begin{enumerate}
\item
$H$ is \textit{spacelike}, i.e. $\langle H,H \rangle >0$.

\item
$H$ is \textit{timelike}, i.e. $\langle H,H \rangle <0$.

\item
$H$ is \textit{lightlike}, i.e. $\langle H,H \rangle =0$.
\end{enumerate}

\vskip 2mm In this section we shall consider spacelike surfaces
whose mean curvature vector at any point is a non-zero spacelike
vector. Let $x,y$ be the principal tangent vector fields. We
denote by $b$ the unit normal vector field
$b=\displaystyle{\frac{H}{\sqrt{\langle H, H \rangle}}}$. We have
$\langle b, b\rangle=1$ and $b$ is collinear with $\sigma(x,x)$
and $\sigma(y,y)$. Denote by $l$ the unit normal vector field such
that $\{x,y,b,l\}$ is a positively oriented orthonormal frame
field in $\R^4_1$ (hence $\langle l, l\rangle=-1$). Thus we obtain
a geometrically determined orthonormal frame field $\{x,y,b,l\}$
at each point $p \in M^2$. With respect to the frame field
$\{x,y,b,l\}$ we have the following formulas:
$$\begin{array}{l}
\vspace{2mm}
\sigma(x,x) = \nu_1\,b; \\
\vspace{2mm}
\sigma(x,y) = \lambda\,b - \mu\,l;  \\
\vspace{2mm} \sigma(y,y) = \nu_2\,b,
\end{array}\leqno{(5.1)}$$
where $\nu_1, \nu_2, \lambda, \mu$ are invariant functions, $\nu_1
= \langle \sigma(x,x), b\rangle$, $\nu_2 = \langle \sigma(y,y),
b\rangle$, $\lambda = \langle \sigma(x,y), b\rangle$, $\mu =
\langle \sigma(x,y), l\rangle$.

The invariants $k$, $\varkappa$, and the Gauss curvature $K$ of
$M^2$ are expressed by the functions $\nu_1, \nu_2, \lambda, \mu$
as follows:
$$k = - 4\nu_1\,\nu_2\,\mu^2, \quad \quad \varkappa = (\nu_1-\nu_2)\mu, \quad \quad
K = \nu_1\,\nu_2- \lambda^2 + \mu^2.\leqno(5.2)$$ Since
$\varkappa^2 - k > 0$, equalities (5.2) imply that $\mu \neq 0$.

The normal mean curvature vector field of $M^2$ is $H =
\ds{\frac{\nu_1 + \nu_2}{2}\, b}$. Taking into account (5.2) we
obtain that the length $\Vert H \Vert$ of the mean curvature
vector field is given by the formula
$$\Vert H \Vert = \displaystyle{\frac{\sqrt{\varkappa^2-k}}{2 |\mu |}},$$
which shows that $|\mu|$ is expressed by the invariants $k$,
$\varkappa$ and the mean curvature function.

Now we shall discuss the geometric meaning of the invariant
$\lambda$. Let $M$ be an $n$-dimensional submanifold of
$(n+m)$-dimensional Riemannian manifold $\widetilde{M}$ and $\xi$
be a normal vector field of $M$. In \cite{Chen1} B.-Y. Chen
defined the \emph{allied vector field} $a(\xi)$ of $\xi$ by the
formula
$$a(\xi) = \ds{\frac{\|\xi\|}{n} \sum_{k=2}^m \{\tr(A_1 A_k)\}\xi_k},$$
where $\{\xi_1 = \ds{\frac{\xi}{\|\xi\|}},\xi_2, \xi_m \}$ is an
orthonormal base of the normal space of $M$, and $A_i = A_{\xi_i},
\,\, i = 1,\dots, m$ is the shape operator with respect to
$\xi_i$. In particular, the allied vector field $a(H)$ of the mean
curvature vector field $H$ is a well-defined normal vector field
which is orthogonal to $H$. It is called the \emph{allied mean
curvature vector field} of $M$ in $\widetilde{M}$. B.-Y. Chen
defined  the $\mathcal{A}$-submanifolds to be those submanifolds
of $\widetilde{M}$ for which
 $a(H)$ vanishes identically \cite{Chen1}.
In \cite{GVV1}, \cite{GVV2} the $\mathcal{A}$-submanifolds are
called \emph{Chen submanifolds}. It is easy to see that minimal
submanifolds, pseudo-umbilical submanifolds and hypersurfaces are
Chen submanifolds. These Chen submanifolds are said to be trivial
$\mathcal{A}$-submanifolds. Now let $M^2$ be a spacelike surface
in $\R^4_1$ with spacelike mean curvature vector field. Applying
the definition of the allied mean curvature vector field from
equalities (5.1) we get
$$a(H) = \ds{\frac{\nu_1 + \nu_2}{2} \,\lambda \mu \, l}= \ds{\frac{\sqrt{\varkappa^2-k}}{2} \, \lambda \,l}.$$
Hence, if $M^2$ is free of minimal points, then $a(H) = 0$ if and
only if $\lambda = 0$. This gives the geometric meaning of the
invariant $\lambda$. It is clear that $M^2$ is a non-trivial  Chen
surface if and only if the invariant $\lambda$ is zero.

\vskip 2mm

With respect to the geometrically determined orthonormal frame
field $\{x,y,b,l\}$ we have the following Frenet-type derivative
formulas of $M^2$:
$$\begin{array}{ll}
\vspace{2mm} \nabla'_xx=\quad \quad \quad \gamma_1\,y+\,\nu_1\,b;
& \qquad
\nabla'_xb=-\nu_1\,x-\lambda\,y\quad\quad \quad -\beta_1\,l;\\
\vspace{2mm} \nabla'_xy=-\gamma_1\,x\quad \quad \; + \, \lambda\,b
\; - \mu\,l;  & \qquad
\nabla'_yb=-\lambda\,x - \; \nu_2\,y\quad\quad \;\;\, -\beta_2\,l;\\
\vspace{2mm} \nabla'_yx=\quad\quad \;\, -\gamma_2\,y  + \lambda\,b
\; -\mu\,l;  & \qquad
\nabla'_xl= \quad \quad \quad \;-\mu\,y-\beta_1\,b;\\
\vspace{2mm} \nabla'_yy=\;\;\gamma_2\,x \quad\quad\;\;\,
+\nu_2\,b; & \qquad \nabla'_yl=-\mu\,x \quad \quad \quad
\;-\beta_2\,b,
\end{array}\leqno{(5.3)}$$
where  $\gamma_1 = - y(\ln \sqrt{E}), \,\, \gamma_2 = - x(\ln
\sqrt{G})$, $\beta_1 = \langle \nabla'_xb, l\rangle$, $\beta_2 =
\langle \nabla'_yb, l\rangle$.

Using that $R'(x,y,x) = 0$, $R'(x,y,y) = 0$, and $R'(x,y,b) = 0$,
from (5.3) we get the following integrability conditions:
$$\begin{array}{l}
\vspace{2mm} \nu_1 \,\nu_2 - \lambda^2 + \mu^2 = x(\gamma_2) +
y(\gamma_1) - \left((\gamma_1)^2
+ (\gamma_2)^2\right);\\
\vspace{2mm}
2\mu\, \gamma_2 + \nu_1\,\beta_2 - \lambda\,\beta_1 = x(\mu);\\
\vspace{2mm}
2\mu\, \gamma_1 - \lambda\,\beta_2 + \nu_2\,\beta_1 = y(\mu);\\
\vspace{2mm} 2\lambda\, \gamma_2 - \mu\,\beta_1 - (\nu_1 -
\nu_2)\,\gamma_1 =
x(\lambda) - y(\nu_1);\\
\vspace{2mm} 2\lambda\, \gamma_1 - \mu\,\beta_2 + (\nu_1 -
\nu_2)\,\gamma_2 =
- x(\nu_2) + y(\lambda);\\
\gamma_1\,\beta_1 - \gamma_2\,\beta_2 + (\nu_1 - \nu_2)\,\mu  = -
x(\beta_2) + y(\beta_1).
\end{array}$$

Having in mind that $x = \displaystyle{\frac{z_u}{\sqrt{E}}, \, y
= \frac{z_v}{\sqrt{G}}}$, we can rewrite the above equalities in
the following way:
$$\begin{array}{l}
\vspace{2mm} \nu_1 \,\nu_2 - \lambda^2 + \mu^2 =
\displaystyle{\frac{1}{\sqrt{E}}\,(\gamma_2)_u
+ \frac{1}{\sqrt{G}}\,(\gamma_1)_v - \left((\gamma_1)^2 + (\gamma_2)^2\right)};\\
\vspace{2mm} 2\mu\, \gamma_2 + \nu_1\,\beta_2 - \lambda\,\beta_1 =
\displaystyle{\frac{1}{\sqrt{E}} \, \mu_u};\\
\vspace{2mm} 2\mu\, \gamma_1 - \lambda\,\beta_2 + \nu_2\,\beta_1 =
\displaystyle{\frac{1}{\sqrt{G}}\, \mu_v};\\
\vspace{2mm} 2\lambda\, \gamma_2 - \mu\,\beta_1 - (\nu_1 -
\nu_2)\,\gamma_1 =
\displaystyle{\frac{1}{\sqrt{E}}\, \lambda_u - \frac{1}{\sqrt{G}}\,(\nu_1)_v};\\
\vspace{2mm} 2\lambda\, \gamma_1 - \mu\,\beta_2 + (\nu_1 -
\nu_2)\,\gamma_2 =
\displaystyle{ - \frac{1}{\sqrt{E}}\,(\nu_2)_u + \frac{1}{\sqrt{G}}\,\lambda_v};\\
\gamma_1\,\beta_1 - \gamma_2\,\beta_2 + (\nu_1 - \nu_2)\,\mu  =
\displaystyle{ - \frac{1}{\sqrt{E}}\,(\beta_2)_u +
\frac{1}{\sqrt{G}}\,(\beta_1)_v}.
\end{array}$$

The condition $\mu_u \,\mu_v \neq 0$ is equivalent to $(2\mu\,
\gamma_2 + \nu_1\,\beta_2 - \lambda\,\beta_1) (2\mu\, \gamma_1 -
\lambda\,\beta_2 + \nu_2\,\beta_1) \neq 0$. So, if $\mu_u \,\mu_v
\neq 0$, then
$$\sqrt{E} = \displaystyle{\frac{\mu_u}
{2\mu\, \gamma_2 + \nu_1\,\beta_2 - \lambda\,\beta_1}}; \quad
\sqrt{G} = \displaystyle{\frac{\mu_v}{2\mu\, \gamma_1 -
\lambda\,\beta_2 + \nu_2\,\beta_1}}.$$

\vskip 3mm We shall prove the following  Bonnet-type theorem for
spacelike surfaces in $\R^4_1$ whose mean curvature vector at any
point is a non-zero spacelike vector.

\begin{thm}\label{T:Main Theorem}
Let $\gamma_1, \, \gamma_2, \, \nu_1,\, \nu_2, \, \lambda, \, \mu,
\, \beta_1, \beta_2$ be smooth functions, defined in a domain
$\mathcal{D}, \,\, \mathcal{D} \subset {\R}^2$, and satisfying the
conditions
$$\begin{array}{l}
\vspace{2mm}
\displaystyle{\frac{\mu_u}{2\mu\, \gamma_2 + \nu_1\,\beta_2 - \lambda\,\beta_1}}>0;\\
\vspace{2mm}
\displaystyle{\frac{\mu_v}{2\mu\, \gamma_1 - \lambda\,\beta_2 + \nu_2\,\beta_1}}>0;\\
\vspace{2mm}
- \gamma_1 \sqrt{E} \sqrt{G} = (\sqrt{E})_v;\\
\vspace{2mm}
- \gamma_2 \sqrt{E} \sqrt{G} = (\sqrt{G})_u;\\
\vspace{2mm} \nu_1 \,\nu_2 - \lambda^2 + \mu^2 =
\displaystyle{\frac{1}{\sqrt{E}}\,(\gamma_2)_u
+ \frac{1}{\sqrt{G}}\,(\gamma_1)_v - \left((\gamma_1)^2 + (\gamma_2)^2\right)};\\
\vspace{2mm} 2\lambda\, \gamma_2 - \mu\,\beta_1 - (\nu_1 -
\nu_2)\,\gamma_1 =
\displaystyle{\frac{1}{\sqrt{E}}\,\lambda_u - \frac{1}{\sqrt{G}}\,(\nu_1)_v};\\
\vspace{2mm} 2\lambda\, \gamma_1 - \mu\,\beta_2 + (\nu_1 -
\nu_2)\,\gamma_2 =
\displaystyle{ - \frac{1}{\sqrt{E}}\,(\nu_2)_u + \frac{1}{\sqrt{G}}\,\lambda_v};\\
\gamma_1\,\beta_1 - \gamma_2\,\beta_2 + (\nu_1 - \nu_2)\,\mu  =
\displaystyle{ - \frac{1}{\sqrt{E}}\,(\beta_2)_u +
\frac{1}{\sqrt{G}}\,(\beta_1)_v},
\end{array}\leqno{(5.4)}$$
where $\sqrt{E} = \displaystyle{\frac{\mu_u}{2\mu\, \gamma_2 +
\nu_1\,\beta_2 - \lambda\,\beta_1}}$, $\sqrt{G} =
\displaystyle{\frac{\mu_v} {2\mu\, \gamma_1 - \lambda\,\beta_2 +
\nu_2\,\beta_1}}$. Let $\{x_0, \, y_0, \, b_0,\, l_0\}$ be an
orthonormal frame at a point $p_0 \in \R^4_1$. Then there exist a
subdomain ${\mathcal{D}}_0 \subset \mathcal{D}$ and a unique
spacelike surface $M^2: z = z(u,v), \,\, (u,v) \in
{\mathcal{D}}_0$, whose mean curvature vector at any point is a
non-zero spacelike vector. Moreover, $M^2$ passes through $p_0$,
the functions
 $\gamma_1, \, \gamma_2, \, \nu_1,\, \nu_2, \,
\lambda, \, \mu, \, \beta_1, \beta_2$ are the geometric functions
of $M^2$ and $\{x_0, \, y_0, \, b_0,\, l_0\}$ is the geometric
frame of $M^2$ at the point $p_0$.
\end{thm}

\vskip 2mm \noindent \emph{Proof:} We consider the following
system of partial differential equations for the unknown vector
functions $x = x(u,v), \, y = y(u,v), \,b = b(u,v), \,l = l(u,v)$
in $\R^4_1$:
$$\begin{array}{ll}
\vspace{2mm} x_u = \sqrt{E}\, \gamma_1\, y + \sqrt{E}\, \nu_1\, b
& \qquad x_v =
- \sqrt{G}\, \gamma_2\, y + \sqrt{G}\, \lambda\, b - \sqrt{G}\, \mu\, l\\
\vspace{2mm} y_u = - \sqrt{E}\, \gamma_1\, x + \sqrt{E}\,
\lambda\, b
- \sqrt{E}\, \mu\, l  & \qquad y_v = \sqrt{G}\, \gamma_2\, x + \sqrt{G}\, \nu_2\, b \\
\vspace{2mm} b_u = - \sqrt{E}\, \nu_1\, x - \sqrt{E}\, \lambda\, y
- \sqrt{E}\, \beta_1\, l  &
\qquad b_v = - \sqrt{G}\, \lambda\, x - \sqrt{G}\, \nu_2\, y - \sqrt{G}\, \beta_2\, l \\
\vspace{2mm} l_u = - \sqrt{E}\, \mu\, y - \sqrt{E}\, \beta_1\, b &
\qquad l_v = - \sqrt{G}\, \mu\, x - \sqrt{G}\, \beta_2\, b
\end{array}\leqno{(5.5)}$$
We denote
$$Z =
\left(%
\begin{array}{c}
  x \\
  y \\
  b \\
  l \\
\end{array}%
\right); \quad
A = \sqrt{E} \left(%
\begin{array}{cccc}
  0 & \gamma_1 & \nu_1 & 0 \\
  -\gamma_1 & 0 & \lambda & - \mu \\
  -\nu_1 & -\lambda & 0 & - \beta_1 \\
  0 & -\mu & -\beta_1 & 0 \\
\end{array}%
\right); \quad B = \sqrt{G}
\left(%
\begin{array}{cccc}
  0 & -\gamma_2 & \lambda & - \mu \\
  \gamma_2 & 0 & \nu_2 & 0 \\
  -\lambda & -\nu_2 & 0 & - \beta_2 \\
  -\mu & 0 & -\beta_2 & 0 \\
\end{array}%
\right).$$ Then system (5.5) can be rewritten in the form:
$$\begin{array}{l}
\vspace{2mm}
Z_u = A\,Z,\\
\vspace{2mm} Z_v = B\,Z.
\end{array}\leqno{(5.6)}$$
The integrability conditions of (5.6) are
$$Z_{uv} = Z_{vu},$$
i.e.
$$\displaystyle{\frac{\partial a_i^k}{\partial v} - \frac{\partial b_i^k}{\partial u}
+ \sum_{j=1}^{4}(a_i^j\,b_j^k - b_i^j\,a_j^k) = 0, \quad i,k = 1,
\dots, 4,} \leqno(5.7)$$ where $a_i^j$ and $b_i^j$ are the
elements of the matrices $A$ and $B$. Using (5.4) we obtain that
 equalities (5.7) are fulfilled. Hence, there exists a subset
$\mathcal{D}_1 \subset \mathcal{D}$ and unique vector functions $x
= x(u,v), \, y = y(u,v), \,b = b(u,v), \,l = l(u,v), \,\, (u,v)
\in \mathcal{D}_1$, which satisfy system (5.5) and the conditions
$$x(u_0,v_0) = x_0, \quad y(u_0,v_0) = y_0, \quad b(u_0,v_0) = b_0, \quad l(u_0,v_0) = l_0.$$

We shall prove that $x(u,v), \, y(u,v), \,b(u,v), \,l(u,v)$ form
an orthonormal frame in $\R^4_1$ for each $(u,v) \in
\mathcal{D}_1$. Let us consider the following functions:
$$\begin{array}{lll}
\vspace{2mm}
  \varphi_1 = \langle x,x \rangle - 1; & \qquad \varphi_5 =
  \langle x,y \rangle; & \qquad \varphi_8 = \langle y,b \rangle; \\
\vspace{2mm}
  \varphi_2 = \langle y, y \rangle - 1; & \qquad \varphi_6 =
  \langle x,b \rangle; & \qquad \varphi_9 = \langle y,l \rangle; \\
\vspace{2mm}
  \varphi_3 = \langle b, b \rangle - 1; & \qquad \varphi_7 =
  \langle x,l \rangle; & \qquad \varphi_{10} = \langle b,l \rangle; \\
\vspace{2mm}
  \varphi_4 = \langle l,l \rangle + 1; &   &  \\
\end{array}$$
defined for each $(u,v) \in \mathcal{D}_1$. Using that $x(u,v), \,
y(u,v), \,b(u,v), \,l(u,v)$ satisfy (5.5), we obtain  the system
$$\begin{array}{lll}
\vspace{2mm}
\displaystyle{\frac{\partial \varphi_i}{\partial u} = \alpha_i^j \, \varphi_j},\\
\vspace{2mm} \displaystyle{\frac{\partial \varphi_i}{\partial v} =
\beta_i^j \, \varphi_j};
\end{array} \qquad i = 1, \dots, 10, \leqno{(5.8)}$$
where $\alpha_i^j, \beta_i^j, \,\, i,j = 1, \dots, 10$ are
functions of $(u,v) \in \mathcal{D}_1$. System (5.8) is a linear
system of partial differential equations for the functions
$\varphi_i(u,v), \,\,i = 1, \dots, 10, \,\,(u,v) \in
\mathcal{D}_1$, satisfying $\varphi_i(u_0,v_0) = 0, \,\,i = 1,
\dots, 10$. Hence, $\varphi_i(u,v) = 0, \,\,i = 1, \dots, 10$ for
each $(u,v) \in \mathcal{D}_1$. Consequently, the vector functions
$x(u,v), \, y(u,v), \,b(u,v), \,l(u,v)$ form an orthonormal frame
in $\R^4_1$ for each $(u,v) \in \mathcal{D}_1$.

Now, let us consider the  system
$$\begin{array}{lll}
\vspace{2mm}
z_u = \sqrt{E}\, x\\
\vspace{2mm} z_v = \sqrt{G}\, y
\end{array}\leqno{(5.9)}$$
of partial differential equations for the vector function
$z(u,v)$. Using (5.4) and (5.5) we get that the integrability
conditions $z_{uv} = z_{vu}$ of system (5.9)
 are fulfilled. Hence,  there exists a subset $\mathcal{D}_0 \subset \mathcal{D}_1$ and
a unique vector function $z = z(u,v)$, defined for $(u,v) \in
\mathcal{D}_0$ and satisfying $z(u_0, v_0) = p_0$.

Consequently, the surface $M^2: z = z(u,v), \,\, (u,v) \in
\mathcal{D}_0$ satisfies the assertion of the theorem. \qed

\section{Spacelike surfaces with timelike mean curvature vector field}\label{S:timelikeH}

Now we shall consider spacelike surfaces with timelike mean
curvature vector field, i.e. $H \neq 0$, $\langle H, H\rangle <
0$. Let $x,y$ be the principal tangent vector fields. We denote by
$l$ the unit normal vector field $l=\displaystyle{-
\frac{H}{\sqrt{- \langle H, H\rangle}}}$. We have $\langle l,
l\rangle= - 1$ and $l$ is collinear with $\sigma(x,x)$ and
$\sigma(y,y)$. Denote by $b$ the unit normal vector field
($\langle b, b\rangle=1$) such that the quadruple $\{x,y,b,l\}$ is
a positively oriented orthonormal frame field in $\R^4_1$. Thus we
obtain a geometrically determined orthonormal frame field
$\{x,y,b,l\}$ at each point $p \in M^2$. With respect to the frame
field $\{x,y,b,l\}$ we have the following formulas:
$$\begin{array}{l}
\vspace{2mm}
\sigma(x,x) = \quad \;\,- \nu_1\,l; \\
\vspace{2mm}
\sigma(x,y) = \mu\,b - \lambda\,l;  \\
\vspace{2mm} \sigma(y,y) = \quad \;\,- \nu_2\,l,
\end{array}$$
where $\nu_1, \nu_2, \lambda, \mu$ are invariant functions, $\nu_1
= \langle \sigma(x,x), l\rangle$, $\nu_2 = \langle \sigma(y,y),
l\rangle$, $\lambda = \langle \sigma(x,y), l\rangle$, $\mu =
\langle \sigma(x,y), b\rangle$.

The invariants $k$, $\varkappa$, and the Gauss curvature $K$ of
$M^2$ are expressed by the functions $\nu_1, \nu_2, \lambda, \mu$
as follows:
$$k = - 4\nu_1\,\nu_2\,\mu^2, \quad \quad \varkappa = (\nu_1-\nu_2)\mu, \quad \quad
K = -\nu_1\,\nu_2 + \lambda^2 - \mu^2.\leqno(6.1)$$ Since
$\varkappa^2 - k > 0$,  equalities (6.1) imply that $\mu \neq 0$.
The normal mean curvature vector field of $M^2$ is $H =  \ds{-
\frac{\nu_1 + \nu_2}{2}\, l}$. The allied mean curvature vector
field is
$$a(H) = \ds{\frac{\nu_1 + \nu_2}{2} \, \lambda \mu\,b}.$$
As in the previous section we see that $M^2$ is a non-trivial Chen
surface if and only if the invariant $\lambda$ is zero.

Now the  Frenet-type derivative formulas of $M^2$ are:
$$\begin{array}{ll}
\vspace{2mm} \nabla'_xx=\quad \quad \quad \gamma_1\,y \qquad \quad
- \nu_1\,l;  & \qquad
\nabla'_xb= \quad \quad \quad -\mu\,y \qquad \;\;-\beta_1\,l;\\
\vspace{2mm} \nabla'_xy=-\gamma_1\,x\quad \quad \; + \; \mu\,b \;
- \lambda\,l;  & \qquad
\nabla'_yb=-\mu\,x \quad \quad \quad \qquad\;\; -\beta_2\,l;\\
\vspace{2mm} \nabla'_yx=\quad\quad \;-\gamma_2\,y \; + \mu\,b \;
-\lambda\,l;  & \qquad
\nabla'_xl=-\nu_1\,x-\lambda\,y -\beta_1\,b;\\
\vspace{2mm} \nabla'_yy=\;\;\gamma_2\,x \quad\quad\quad\qquad \;\;
-\nu_2\,l; & \qquad \nabla'_yl=-\lambda\,x - \;
\nu_2\,y-\beta_2\,b,
\end{array}\leqno{(6.2)}$$
where  $\gamma_1 = - y(\ln \sqrt{E}), \,\, \gamma_2 = - x(\ln
\sqrt{G})$, $\beta_1 = \langle \nabla'_xb, l\rangle$, $\beta_2 =
\langle \nabla'_yb, l\rangle$.

Using that $R'(x,y,x) = 0$, $R'(x,y,y) = 0$, and $R'(x,y,b) = 0$,
from (6.2) we get the following integrability conditions:
$$\begin{array}{l}
\vspace{2mm} -\nu_1 \,\nu_2 + \lambda^2 - \mu^2 = x(\gamma_2) +
y(\gamma_1) - \left((\gamma_1)^2
+ (\gamma_2)^2\right);\\
\vspace{2mm}
2\mu\, \gamma_2 + \nu_1\,\beta_2 - \lambda\,\beta_1 = x(\mu);\\
\vspace{2mm}
2\mu\, \gamma_1 - \lambda\,\beta_2 + \nu_2\,\beta_1 = y(\mu);\\
\vspace{2mm} 2\lambda\, \gamma_2 - \mu\,\beta_1 - (\nu_1 -
\nu_2)\,\gamma_1 =
x(\lambda) - y(\nu_1);\\
\vspace{2mm} 2\lambda\, \gamma_1 - \mu\,\beta_2 + (\nu_1 -
\nu_2)\,\gamma_2 =
- x(\nu_2) + y(\lambda);\\
\gamma_1\,\beta_1 - \gamma_2\,\beta_2 - (\nu_1 - \nu_2)\,\mu  = -
x(\beta_2) + y(\beta_1).
\end{array}$$

Again the condition $\mu_u \,\mu_v \neq 0$ is equivalent to
$(2\mu\, \gamma_2 + \nu_1\,\beta_2 - \lambda\,\beta_1) (2\mu\,
\gamma_1 - \lambda\,\beta_2 + \nu_2\,\beta_1) \neq 0$. So, if
$\mu_u \,\mu_v \neq 0$, then $\sqrt{E} =
\displaystyle{\frac{\mu_u} {2\mu\, \gamma_2 + \nu_1\,\beta_2 -
\lambda\,\beta_1}}; \,\, \sqrt{G} =
\displaystyle{\frac{\mu_v}{2\mu\, \gamma_1 - \lambda\,\beta_2 +
\nu_2\,\beta_1}}.$ \vskip 3mm

In a similar way as in Section \ref{S:spacelikeH} we prove the
following Bonnet-type theorem for spacelike surfaces in $\R^4_1$
with timelike mean curvature vector field.

\begin{thm}\label{T:Main Theorem-2}
Let $\gamma_1, \, \gamma_2, \, \nu_1,\, \nu_2, \, \lambda, \, \mu,
\, \beta_1, \beta_2$ be smooth functions, defined in a domain
$\mathcal{D}, \,\, \mathcal{D} \subset {\R}^2$, and satisfying the
conditions
$$\begin{array}{l}
\vspace{2mm}
\displaystyle{\frac{\mu_u}{2\mu\, \gamma_2 + \nu_1\,\beta_2 - \lambda\,\beta_1}}>0;\\
\vspace{2mm}
\displaystyle{\frac{\mu_v}{2\mu\, \gamma_1 - \lambda\,\beta_2 + \nu_2\,\beta_1}}>0;\\
\vspace{2mm}
- \gamma_1 \sqrt{E} \sqrt{G} = (\sqrt{E})_v;\\
\vspace{2mm}
- \gamma_2 \sqrt{E} \sqrt{G} = (\sqrt{G})_u;\\
\vspace{2mm} - \nu_1 \,\nu_2 + \lambda^2 - \mu^2 =
\displaystyle{\frac{1}{\sqrt{E}}\,(\gamma_2)_u
+ \frac{1}{\sqrt{G}}\,(\gamma_1)_v - \left((\gamma_1)^2 + (\gamma_2)^2\right)};\\
\vspace{2mm} 2\lambda\, \gamma_2 - \mu\,\beta_1 - (\nu_1 -
\nu_2)\,\gamma_1 =
\displaystyle{\frac{1}{\sqrt{E}}\,\lambda_u - \frac{1}{\sqrt{G}}\,(\nu_1)_v};\\
\vspace{2mm} 2\lambda\, \gamma_1 - \mu\,\beta_2 + (\nu_1 -
\nu_2)\,\gamma_2 =
\displaystyle{ - \frac{1}{\sqrt{E}}\,(\nu_2)_u + \frac{1}{\sqrt{G}}\,\lambda_v};\\
\gamma_1\,\beta_1 - \gamma_2\,\beta_2 - (\nu_1 - \nu_2)\,\mu  =
\displaystyle{ - \frac{1}{\sqrt{E}}\,(\beta_2)_u +
\frac{1}{\sqrt{G}}\,(\beta_1)_v},
\end{array}$$
where $\sqrt{E} = \displaystyle{\frac{\mu_u}{2\mu\, \gamma_2 +
\nu_1\,\beta_2 - \lambda\,\beta_1}}$, $\sqrt{G} =
\displaystyle{\frac{\mu_v} {2\mu\, \gamma_1 - \lambda\,\beta_2 +
\nu_2\,\beta_1}}$. Let $\{x_0, \, y_0, \, b_0,\, l_0\}$ be an
orthonormal frame at a point $p_0 \in \R^4_1$. Then there exist a
subdomain ${\mathcal{D}}_0 \subset \mathcal{D}$ and a unique
spacelike surface $M^2: z = z(u,v), \,\, (u,v) \in
{\mathcal{D}}_0$, passing through $p_0$,  with timelike mean
curvature vector field, such that $\gamma_1, \, \gamma_2, \,
\nu_1,\, \nu_2, \, \lambda, \, \mu, \, \beta_1, \beta_2$ are the
geometric functions of $M^2$ and $\{x_0, \, y_0, \, b_0,\, l_0\}$
is the geometric frame of $M^2$ at the point $p_0$.
\end{thm}

\vskip 5mm
\section{Examples}\label{S:Examples}

In this section we shall apply our theory to a special class of
spacelike surfaces in $\R^4_1$. In  \cite{M} C. Moore studied
general rotational surfaces in $\R^4$. In \cite{GM2,GM4} we
considered a special case of such surfaces, given by
$$\mathcal{M}: z(u,v) = \left( f(u) \cos\alpha v, f(u) \sin \alpha v, g(u) \cos \beta v, g(u) \sin \beta v \right), \leqno{(7.1)}$$
where $u \in J \subset \R, \,\,  v \in [0; 2\pi)$, $f(u)$ and
$g(u)$ are smooth functions, satisfying $\alpha^2 f^2+ \beta^2 g^2
> 0 , \,\, f'\,^2+ g'\,^2 > 0$, and $\alpha, \beta$ are positive
constants. These surfaces are general rotational surfaces in the
sense of C. Moore with plane meridian curves. Here we shall
consider a class of  spacelike surfaces in $\R^4_1$ which are
analogous to (7.1).

\vskip 2mm \noindent \textbf{Example 1.} Let us consider the
surface $\mathcal{M}_1$ parameterized by
$$\mathcal{M}_1: z(u,v) = \left( f(u) \cos \alpha v, f(u) \sin \alpha v, g(u) \cosh \beta v, g(u) \sinh \beta v \right),$$
where $f(u)$ and $g(u)$ are smooth functions, satisfying $\alpha^2
f^2(u)- \beta^2 g^2(u) > 0$, $f'\,^2(u)+ g'\,^2(u) > 0$, $u \in J
\subset \R$ and $\alpha, \beta$ are positive constants;  $v \in
[0; 2\pi)$. The tangent space of $\mathcal{M}_1$ is spanned by the
vector fields
$$\begin{array}{l}
\vspace{2mm}
z_u = \left(f'(u) \cos \alpha v, f'(u) \sin \alpha v, g'(u) \cosh \beta v, g'(u) \sinh \beta v \right),\\
\vspace{2mm} z_v = \left(- \alpha f(u) \sin \alpha v, \alpha f(u)
\cos \alpha v, \beta g(u) \sinh \beta v, \beta g(u) \cosh \beta v
\right).
\end{array}$$
The coefficients of the first fundamental form of $\mathcal{M}_1$
are
$$E = f'\,^2(u)+ g'\,^2(u); \qquad F = 0; \qquad G =\alpha^2 f^2(u)- \beta^2
g^2(u).$$ $\mathcal{M}_1$ is a spacelike surface in $\R^4_1$. We
choose the following normal frame field of $\mathcal{M}_1$:
$$\begin{array}{l}
\vspace{2mm}
n_1 = \ds{\frac{1}{\sqrt{f'\,^2 + g'\,^2}}\left(g' \cos \alpha v, g' \sin \alpha v, - f' \cosh \beta v, - f' \sinh \beta v \right)};\\
\vspace{2mm} n_2 = \ds{\frac{1}{\sqrt{\alpha^2 f^2 - \beta^2
g^2}}\left( - \beta g \sin \alpha v, \beta g \cos \alpha v, \alpha
f \sinh \beta v,  \alpha f \cosh \beta v \right)}.
\end{array}$$
We have $\langle n_1, n_1 \rangle = 1,\,\, \langle n_2, n_2
\rangle = -1$. Calculating the second partial derivatives of
$z(u,v)$ we find the functions $c_{ij}^k$ and get the coefficients
$L$, $M$, $N$ of the second fundamental form of $\mathcal{M}_1$:
$$L = \ds{\frac{2 \alpha \beta (g f' - f g') (g' f'' - f' g'')}{(f'\,^2 + g'\,^2)(\alpha^2 f^2 - \beta^2 g^2)}}; \qquad M = 0; \qquad
N = \ds{\frac{2\alpha \beta (g f' - f g') (\alpha^2 fg' + \beta^2
gf')}{(f'\,^2 + g'\,^2)(\alpha^2 f^2 - \beta^2 g^2)}}.$$
Consequently, the invariants $k$, $\varkappa$, and $K$ of
$\mathcal{M}_1$ are expressed as follows:
$$k = \ds{\frac{4 \alpha^2 \beta^2 (g f' - f g')^2 (g' f'' - f' g'') (\alpha^2 f g' + \beta^2 g f')}{(f'\,^2 + g'\,^2)^3 (\alpha^2 f^2 - \beta^2 g^2)^3 }};$$

$$\varkappa =  \ds{\frac{\alpha \beta (g f' - f g')}{(f'\,^2 + g'\,^2)^2 (\alpha^2 f^2 - \beta^2 g^2)^2} \,
[(\alpha^2 f^2 - \beta^2  g^2)(g' f'' - f' g'') + (f'\,^2 +
g'\,^2) ( \alpha^2 f g' + \beta^2 g f') ]};$$

$$K =  \ds{\frac{- (\alpha^2 f^2 - \beta^2  g^2)(\alpha^2 f g' + \beta^2 g f')(g' f'' - f' g'') + \alpha^2 \beta^2 (f'\,^2 + g'\,^2) (g f' - f g')^2}{(f'\,^2 + g'\,^2)^2 (\alpha^2 f^2 - \beta^2 g^2)^2 } \,}.$$

The mean curvature vector field  $H$ is collinear with $n_1$, and
hence $\mathcal{M}_1$ is a spacelike surface whose mean curvature
vector at any point is a non-zero spacelike vector. Note that
$\mathcal{M}_1$ is parameterized by principal parameters $(u,v)$.
Denoting $x = \ds{\frac{z_u}{\sqrt{f'\,^2(u)+ g'\,^2(u)}}}, \,\, y
= \ds{\frac{z_v} {\sqrt{\alpha^2 f^2(u)- \beta^2 g^2(u)}}}$ we
obtain the geometric invariant functions in the  Frenet-type
derivative formulas of $\mathcal{M}_1$:
$$\begin{array}{ll}
\vspace{2mm}
\gamma_1 = 0; & \qquad \gamma_2 = \ds{- \frac{\alpha^2 f f' - \beta^2 g g'}{\sqrt{f'\,^2 + g'\,^2}(\alpha^2 f^2 - \beta^2 g^2)}};\\
\vspace{2mm} \nu_1 = \ds{\frac{g' f'' - f' g''}{(f'\,^2 +
g'\,^2)^{\frac{3}{2}}}}; & \qquad
\nu_2 = \ds{- \frac{\alpha^2 f g' + \beta^2 g f'}{\sqrt{f'\,^2 + g'\,^2}(\alpha^2 f^2 - \beta^2 g^2)}};\\
\vspace{2mm}
\lambda = 0; & \qquad  \mu = \ds{\frac{\alpha \beta (g f' - f g')}{\sqrt{f'\,^2 + g'\,^2}(\alpha^2 f^2 - \beta^2 g^2)}};\\
\vspace{2mm}
 \beta_1 = 0; & \qquad  \beta_2 = \ds{\frac{\alpha \beta (f f' + g g')}{\sqrt{f'\,^2 + g'\,^2}(\alpha^2 f^2 - \beta^2 g^2) }}.
\end{array}$$

Since the invariant $\lambda$ is zero, the general rotational
surface $\mathcal{M}_1$ is a  Chen surface. In \cite{Houh} C. Houh
considered a more general class of surfaces of rotational type in
$\R^4_1$ and found a subclass of Chen surfaces.

In the special case when $f(u) = \cos u,\,\, g(u) = \sin u$,
$\alpha=\beta=1$ we obtain a spacelike surface lying on  De Sitter
space $S^3_1 = \{x \in \R^4_1; \langle x,x \rangle = 1\}$ with
invariants
$$k = \ds{-\frac{4}{\cos^2 2u}}; \qquad \varkappa = 0; \qquad K = \ds{\frac{\cos^2 2u + 1}{\cos^2 2u}}.$$
This is an example of a spacelike surface with flat normal
connection and spacelike mean curvature vector field.

\vskip 4mm \noindent \textbf{Example 2.} Now we shall consider the
surface $\mathcal{M}_2$ parameterized by
$$\mathcal{M}_2: z(u,v) = \left( f(u) \cos \alpha v, f(u) \sin \alpha v, g(u) \sinh \beta v, g(u) \cosh \beta v \right),$$
where $f(u)$ and $g(u)$ are smooth functions, satisfying
$f'\,^2(u)- g'\,^2(u) > 0$, $\alpha^2 f^2(u)+ \beta^2 g^2(u) > 0$,
$u \in J \subset \R$ and $\alpha, \beta$ are positive constants;
$v \in [0; 2\pi)$. The coefficients of the first fundamental form
of $\mathcal{M}_2$ are
$$E = f'\,^2(u)- g'\,^2(u); \qquad F = 0; \qquad G =\alpha^2 f^2(u)+ \beta^2 g^2(u),$$
hence $\mathcal{M}_2$ is a spacelike surface in $\R^4_1$. We
choose the following normal frame field of $\mathcal{M}_2$:
$$\begin{array}{l}
\vspace{2mm} n_1 = \ds{\frac{1}{\sqrt{\alpha^2 f^2 + \beta^2
g^2}}\left( \beta g \sin \alpha v,- \beta g \cos \alpha v, \alpha
f \cosh \beta v,  \alpha f \sinh \beta v\right)};\\
\vspace{2mm} n_2 = \ds{\frac{1}{\sqrt{f'\,^2 - g'\,^2}}\left(g'
\cos \alpha v, g' \sin \alpha v,  f' \sinh \beta v, f' \cosh \beta
v  \right)}.
\end{array}$$
We have $\langle n_1, n_1 \rangle = 1,\,\, \langle n_2, n_2
\rangle = -1$. Calculating the coefficients $L$, $M$, $N$ of the
second fundamental form we obtain that the invariants
 $k$, $\varkappa$, and $K$ of $\mathcal{M}_2$ are expressed by the functions $f(u)$, $g(u)$ and their derivatives as
follows:
$$k = \ds{\frac{4 \alpha^2 \beta^2 (g f' - f g')^2 (g' f'' - f' g'') (\alpha^2 f g' + \beta^2 g f')}{(f'\,^2 - g'\,^2)^3 (\alpha^2 f^2 + \beta^2 g^2)^3 }};$$

$$\varkappa =  \ds{\frac{\alpha \beta (g f' - f g')}{(f'\,^2 - g'\,^2)^2 (\alpha^2 f^2 + \beta^2 g^2)^2} \,
[(\alpha^2 f^2 + \beta^2  g^2)(g' f'' - f' g'') + (f'\,^2 -
g'\,^2) ( \alpha^2 f g' + \beta^2 g f') ]};$$

$$K =  \ds{\frac{(\alpha^2 f^2 + \beta^2  g^2)(\alpha^2 f g' + \beta^2 g f')(g' f'' - f' g'') - \alpha^2 \beta^2 (f'\,^2 - g'\,^2) (g f' - f g')^2}{(f'\,^2 - g'\,^2)^2 (\alpha^2 f^2 + \beta^2 g^2)^2 } \,}.$$

In this example the mean curvature vector field  $H$ is collinear
with $n_2$, and hence $\mathcal{M}_2$ is a spacelike surface with
timelike mean curvature vector field. We note that $\mathcal{M}_2$
is parameterized by principal parameters $(u,v)$. The geometric
invariant functions in the  Frenet-type derivative formulas of
$\mathcal{M}_2$ are given below:
$$\begin{array}{ll}
\vspace{2mm}
\gamma_1 = 0; & \qquad \gamma_2 = \ds{- \frac{\alpha^2 f f' + \beta^2 g g'}{\sqrt{f'\,^2 - g'\,^2}(\alpha^2 f^2 + \beta^2 g^2)}};\\
\vspace{2mm} \nu_1 = \ds{\frac{g' f'' - f' g''}{(f'\,^2 -
g'\,^2)^{\frac{3}{2}}}}; & \qquad
\nu_2 = \ds{- \frac{\alpha^2 f g' + \beta^2 g f'}{\sqrt{f'\,^2 - g'\,^2}(\alpha^2 f^2 + \beta^2 g^2)}};\\
\vspace{2mm}
\lambda = 0; & \qquad  \mu = \ds{\frac{\alpha \beta (f g' - g f')}{\sqrt{f'\,^2 - g'\,^2}(\alpha^2 f^2 + \beta^2 g^2)}};\\
\vspace{2mm}
 \beta_1 = 0; & \qquad  \beta_2 = \ds{\frac{\alpha \beta (g g' - f f')}{\sqrt{f'\,^2 - g'\,^2}(\alpha^2 f^2 + \beta^2 g^2) }}.
\end{array}$$
The surface $\mathcal{M}_2$ is a spacelike  Chen surface, since
the invariant $\lambda$ is zero.

If we choose  $f(u) = \sinh u,\,\, g(u) = \cosh u$,
$\alpha=\beta=1$ we obtain a spacelike surface lying on the unit
hyperbolic sphere $H^3_1 = \{x \in \R^4_1; \langle x,x \rangle = -
1\}$ with invariants
$$k = \ds{-\frac{4}{\cosh^2 2u}}; \qquad \varkappa = 0; \qquad K = \ds{- \frac{\cosh^2 2u + 1}{\cosh^2 2u}}.$$
This is  a spacelike surface with flat normal connection and
timelike mean curvature vector field.

\vskip 2mm \textbf{Acknowledgements:} The authors would like to
express their thanks to the referee for his valuable comments and
suggestions. The second author is partially supported by "L.
Karavelov" Civil Engineering Higher School, Sofia, Bulgaria under
Contract No 10/2010.

\end{document}